\documentclass[10pt]{amsart}
\usepackage{latexsym}   
\usepackage{amssymb}    
\usepackage{amsmath}    
\usepackage{amsthm}     
\usepackage[all]{xy}
\usepackage{hyperref}
\newcommand{\pa}{\partial}\newcommand{\al}{\alpha}
\newcommand{\be}{\beta}

\newcommand{\del}{\delta}
\newcommand{\ep}{\epsilon}

\newcommand{\om}{\omega}
\newcommand{\Om}{\Omega}
\newcommand{\ti}{\tilde}

\renewcommand{\thefootnote}

\newtheorem{theorem}{Theorem}[section]

\theoremstyle{definition}

\theoremstyle{remark}

\numberwithin{equation}{section}

\title[On Bianchi's B\"{a}cklund transformation of quadrics]
{On Bianchi's B\"{a}cklund transformation of quadrics}

\author[  Ion I. Dinc\u{a}]{Ion I. Dinc\u{a}}

\thanks{Supported by the University of Bucharest}
\address{Faculty of Mathematics and Informatics,
University of Bucharest,  14 Academiei Str., 010014, Bucharest,
Romania}
 \email{dinca@gta.math.unibuc.ro}

\begin{document}


\begin{abstract}
We investigate basic features of Bianchi's B\"{a}cklund
transformation of quadrics to see if it can be obtained under
weaker assumptions and if it can be generalized to deformations of
other surfaces.
\end{abstract}

\maketitle

\tableofcontents \pagenumbering{arabic}

\section{Introduction}
It has been an open question ever since the classical times to
find all surfaces for which an interesting theory of deformation
(procedures of generating more or less explicit examples depending
on arbitrarily many constants; for example {\it B\"{a}cklund} (B)
transformation) can be build, just as for quadrics (see \S\ 5.17
of R. Calapso's introduction to Vol {\bf 4}, part 1 of Bianchi's
{\it Opere}).

Based on Bianchi's \cite{B3} (also appearing as (\cite{B5},Vol
{\bf 4},(143))) characterization of the (isotropic) singular B
transformation of quadrics with auxiliary surface plane, we
consider basic features of Bianchi's {\it B\"{a}cklund} (B)
transformation of quadrics to see if it can be obtained under
weaker assumptions and if it can be generalized to {\it
deformations} (isometric deformations) of other (classes of)
surfaces.

Due to the frequent use of certain keywords we shall define and
use abbreviated notations.

All computations are local and assumed to be valid on their open
domain of validity without further details; all functions have the
assumed order of differentiability and are assumed to be
invertible, non-zero, etc when required (for all practical
purposes we can assume all functions to be analytic).

By necessity (all quadrics in the real $3$-dimensional Euclidean
space are equivalent from a complex projective point of view) we
shall consider the complexification

$$(\mathbb{C}^3,<.,.>),\ <x,y>:=x^Ty,\ |x|^2:=x^Tx,\ x,y\in\mathbb{C}^3$$
of the real $3$-dimensional Euclidean space.

{\it Isotropic} (null) vectors are those vectors of length $0$;
since most vectors are not isotropic we call a vector simply
vector and we shall emphasize isotropic for isotropic vectors. The
same denomination will apply in other settings: for example we
call quadric a non-degenerate quadric (a quadric projectively
equivalent to the complex unit sphere).

We call {\it surface} any sub-manifold of $\mathbb{C}^3\simeq
\mathbb{R}^6$ of real dimension $2,3$ or $4$ such that all its
{\it complexified tangent spaces} (called henceforth {\it tangent
spaces}) have complex dimension $2$. In this case the distribution
$T\cap iT$ formed by intersecting real tangent spaces $T$ with
$iT$ has constant real dimension (which must respectively be $0,\
2$ or $4$). This distribution is integrable (by its own
definition); by an application of Frobenius one obtains leaves; on
leaves the Nijenhuis condition (with the almost complex structure
being induced by the surrounding space) holds and by
Nirenberg-Newlander this is equivalent (locally) to prescribing
$x:D\rightarrow\mathbb{C}^3$, where $D$ is a domain of
$\mathbb{R}^2$ or $\mathbb{C}\times\mathbb{R}$ or $\mathbb{C}^2$
such that $dx\times\wedge dx\neq 0$ (the $\times$ applies to the
vector structure and the $\wedge$ to the exterior form structure,
so the order does not matter: $\times\wedge=\wedge\times$); from a
practical point of view we need to consider only surfaces that are
real $2$- or $4$-dimensional.

For any two curves there is a developable circumscribed to them;
{\it isotropic developables} are surfaces circumscribed to a curve
and to the circle $C(\infty)$ at $\infty$ (Cayley's absolute);
they are the only surfaces with degenerate $2$-dimensional linear
element. Thus we assume a surface not to be isotropic developable
unless otherwise stated; note however that the importance of
isotropic developables cannot be cast aside, as most isotropic
developables circumscribed to (isotropic) (singular) conics (the
singular part of a doubly ruled (isotropic) plane image of the
unit sphere or of the equilateral paraboloid under an affine
transformation of $\mathbb{C}^3$ with $1$-dimensional kernel
(which is different from the axis of the equilateral paraboloid in
this case)) generate the confocal family of quadrics one of whose
{\it (isotropic) singular quadrics} (including in a general sense
(isotropic) plane of certain pencils of (isotropic) planes)
contains as a singular set the given conic: any quadric of the
confocal family cuts the conic at $4$ points (B\'{e}zout) and the
isotropic rulings of the isotropic developable passing through
those points are all the isotropic rulings of the quadric; three
rulings of a ruling family of a quadric are enough to determine
the quadric (although there are quadrics with less than three
isotropic rulings in both of the ruling families, by inspection
all quadrics except (pseudo-)spheres are uniquely determined by
their finitely many isotropic rulings).

If we let the spectral parameter $z$ vary  in the family of
quadrics confocal to a given one, then we get an affine
correspondence between confocal quadrics (called the {\it Ivory
affinity}) with good metric properties.

In particular it preserves lengths of rulings, so it takes
isotropic rulings and umbilics (their finite intersections; the
remaining $4$ points of intersection of isotropic rulings are
situated on $C(\infty)$ and for either of them having multiplicity
we have quadrics of revolution or {\it Darboux quadrics}) to
isotropic rulings and umbilics; thus as $z$ varies an umbilic
describes a singular conic of a(n isotropic) singular confocal
quadric and an isotropic ruling for each ruling family describes
the isotropic developable circumscribed to the singular conic.

Standard geometric formulae for surfaces in $\mathbb{R}^3$ remain
valid (with their usual denomination) for all surfaces and outside
the locus of isotropic normal directions. Since the {\it
Gau\ss-Weingarten} (GW) and the {\it
Gau\ss-Codazzi-Mainardi-Peterson} (G-CMP) equations for a surface
$x\subset\mathbb{R}^3$ assume only the non-degeneracy of the
linear element of $x$ (and thus the existence of an orthonormal
normal), they are still valid and sufficient to describe the
geometry of surfaces in $\mathbb{C}^3$ almost everywhere; for
example the Lorentz space of signature $(2,1)$ can be realized as
the {\it totally real subspace} $\mathbb{R}^2\times
i\mathbb{R}\subset\mathbb{C}^3$ (the induced scalar product is
real (valued) and non-degenerate).

By the Gau\ss-Bonnet-Peterson fundamental theorem of surfaces the
GW and G-CMP equations suffice to describe the geometry of all
surfaces almost everywhere and surfaces such defined are the
natural completion of the usual real surfaces, although the usual
coordinates used in geometry (asymptotic, conjugate systems,
orthogonal, principal, etc) may not be the ones which clearly
split into purely real and purely complex (for example asymptotic
coordinates on a real surface of positive Gau\ss\ curvature). The
change of coordinates for real $3$-dimensional surfaces still
holds outside bi-holomorphic change of the complex coordinate and
diffeomorphism of the real one, but if the real and complex
coordinates are mixed, they lose the character of being purely
real or complex and the new parametrization may not be holomorphic
in some variables. For example if we have a surface
$x=x(z,t)\subset\mathbb{C}^3,\ (z,t)\in
D\subseteq\mathbb{C}\times\mathbb{R}$ and consider two complex
functions $u=u(z,t),\ v=v(z,t)$ of $z$ and $t$ with
$J:=\frac{du\wedge dv}{dz\wedge dt}=
\begin{vmatrix}\frac{\pa u}{\pa z}&\frac{\pa u}{\pa t}\\
\frac{\pa v}{\pa z}&\frac{\pa v}{\pa t}\end{vmatrix}\neq 0$, then
one can invert $z=z(u,v),\ t=t(u,v)$ at least formally and at the
infinitesimal level by the usual calculus rule of taking the
inverse of the Jacobian: $dz=\frac{\pa z}{\pa u}du+\frac{\pa
z}{\pa v}dv,\ dt=\frac{\pa t}{\pa u}du+\frac{\pa t}{\pa v}dv$.
Although the holomorphic dependence on a variable is lost, the
$2$-dimensionality character remains clear. After using formal
coordinates of geometric meaning to simplify certain computations,
one can return to the purely real and complex coordinates, just as
the use of $z,\bar z$ parametrization on real surfaces. Note that
although on such surfaces the vector fields $\pa_z,\ \pa_{\bar z}$
do not admit integral curves (lines of coordinates), statements
about infinitesimal behavior of such lines remain valid, so one
can assume such lines to exist and derive corresponding results;
thus for all practical purposes we can assume that coordinates
descend upon lines of coordinates on the surface.

Consider Lie's viewpoint: one can replace a surface
$x\subset\mathbb{C}^3$ with a $2$-dimensional distribution of {\it
facets} (pairs of points and planes passing through those points):
the collection of its tangent planes (with the points of tangency
highlighted); thus a facet is the infinitesimal version of a
surface (the integral element $(x,dx)|_{\mathrm{pt}}$ of the
surface). Conversely, a $2$-dimensional distribution of facets is
not always the collection of the tangent planes of a surface (with
the points of tangency highlighted), but the condition that a
$2$-dimensional distribution of facets is integrable (that is it
is the collection of the tangent planes of a leaf (sub-manifold))
does not distinguish between the cases when this sub-manifold is a
surface, curve or point, thus allowing the collapsing of the leaf.

A $3$-dimensional distribution of facets is integrable if it is
the collection of the tangent planes of an $1$-dimensional family
of leaves.

Two {\it rollable} (applicable) surfaces can be {\it rolled}
(applied) one onto the other such that at any instant they meet
tangentially and with the same differential at the tangency point.
The rolling $(x,dx)=(R,t)(x_0,dx_0):=(Rx_0+t,Rdx_0),\
(R,t)\subset\mathbf{O}_3(\mathbb{C})\ltimes\mathbb{C}^3$ of two
{\it applicable} surfaces $x_0,\ x$ (that is $|dx_0|^2=|dx|^2$)
introduces the flat connection form (it encodes the difference of
the second fundamental forms of $x_0,x$ and it being flat encodes
the difference of the G-CMP equations of $x_0,x$).

In rolling various geometric objects (points, curves, surfaces,
facets, etc or families of these) can be rigidly attached to each
point of the rolling surface $x_0$, thus producing rolling objects
(rolling congruences ($2$-dimensional families of objects,
presumed lines unless otherwise stated), rolling distributions of
facets, etc). The study of the rolling problem for surfaces as it
was carried out by the classical geometers (see references
Bianchi, Darboux, Eisenhart) boils down mainly to finding the
features of these rolling objects independent of the shape of the
surface of rolling $x$ (for example envelopes of a rolling
congruence of curves or of a sphere congruence centered on $x_0$,
{\it integrable rolling distributions of facets} (IRDF), {\it
cyclic systems} ($3$-dimensional IRDF such that facets are
centered on circles in the tangent planes of $x_0$ and normal to
the circles)). Note that due to the indeterminacy (in the normal
bundle) of the rolling (one can roll a surface on both sides of
one of its deformations), these rolling objects come in two
families that reflect in the tangent bundle of the rolling
surface.

Conversely, existence (integrability) conditions for such an
object for a particular deformation $x_0$ depend usually on the
linear element of $x_0$ and linearly on the second fundamental
form of $x_0$ (the terms of the second fundamental form appearing
quadratically in the existence conditions become, via Gau\ss\ 's
theorem, dependent of the linear element of $x_0$); to formulate
the conditions for the existence of the rolled objects it is
sufficient to individuate for cancelling the coefficients of the
linearly appearing second fundamental form (note that in the
reflected object these coefficients have changed signs).

Since by infinitesimal rolling in an arbitrary tangential
infinitesimal direction $\del$ an initial facet $\mathcal{F}$
which is common tangent plane to two applicable surfaces is
replaced with an infinitesimally close facet $\mathcal{F}'$ having
in common with $\mathcal{F}$ the direction $\del$, in the actual
rolling problem we have facets centered on each other (the
symmetric {\it tangency configuration} (TC)) and facets centered
on another one $\mathcal{F}$ reflect in $\mathcal{F}$; note that
we assumed a finite law of a general nature as a consequence of an
infinitesimal law via discretization (the converse is clear); see
Bianchi \cite{B3}.

Thus for a theory of deformation of surfaces with the assumption
above we are led to consider, via rolling, certain $4$-dimensional
distributions of facets centered on the tangent planes of the
considered surface $x_0$ and passing through the origin of the
tangent planes (each point of each tangent plane is the center of
finitely many facets) and their rolling counterparts on the
applicable surface $x$.

After a thorough study of infinitesimal laws (in particular
infinitesimal deformations) and their iterations the classical
geometers were led to consider the B transformation (a finite law
of a general nature) as a consequence of an infinitesimal law via
discretization.

The B transformation in the deformation problem naturally appears
by splitting the $4$-dimensional distribution of facets into an
$1$-dimensional family of $3$-dimensional IRDF, thus introducing a
spectral parameter $z\ (B_z)$ (each $3$-dimensional distribution
of facets is integrable (with {\it leaves} $x^1$) regardless of
the shape of the {\it seed} surface $x^0$; note that there are
$4$-dimensional IRDF that do not split into an $1$-dimensional
family of $3$-dimensional IRDF: for example to each point of $x_0$
we associate the tangent planes of a surface).

By imposing the initial collapsing ansatz of leaves of the
$4$-dimensional IRDF to be a $2$-dimensional family of curves
(this ansatz individuates the {\it defining surface} $x_0$) the
$2$-dimensional family of curves naturally splits into an
$1$-dimensional family of {\it auxiliary surfaces} $x_z$
($1$-dimensional family of curves) such that the facets centered
on these auxiliary surfaces form a $3$-dimensional IRDF whose
integrability is independent of the shape of the seed surface
$x^0$ applicable to the surface $x_0^0\subset x_0$. For $x_0$
quadric the auxiliary surfaces $x_z$ are ((isotropic) singular)
confocal quadrics doubly ruled by collapsed leaves and with $z$
being the spectral parameter of the confocal family (this includes
in a general sense the (isotropic) planes of the pencil of
(isotropic) planes through the axis of revolution for $x_0$
quadric of revolution (excluding (pseudo-)spheres) (in this case
the actual isotropic singular quadric of the confocal family is
the two isotropic planes of this pencil) or (isotropic) planes of
the pencil of (isotropic) planes through an isotropic line (which
is the actual isotropic singular quadric of its confocal family)
for $x_0$ Darboux quadric with tangency of order $3$ with
$C(\infty)$).

The case of the singular quadric of the family of confocal
(pseudo-)spheres (the isotropic cone) falls into a different
pattern: the only real cyclic systems with the symmetry of the TC
are obtained for $x_0$ having {\it constant Gau\ss\ curvature}
(CGC) (up to homothety) $-1$, when the leaves also have the same
curvature (Ribaucour, 1870); the transformation from $x_0$ to one
of the leaves is Bianchi's complementary transformation introduced
by him in his Ph. D. thesis in 1879 (thus Bianchi is credited with
the idea of using transformations in the study of (deformations)
of surfaces).

In 1883 B\"{a}cklund introduced his original B transformation of
surfaces of CGC $-1$ by allowing the facets (still centered on
circles in and centered at the center of tangent planes of the
seed) to have constant inclination to the tangent planes of the
seed (thus introducing the spectral parameter).

Lie studied both Bianchi's complementary transformation and
B\"{a}cklund's transformation as arrangements of facets; he proved
their inversion (thus the B transformation can be iterated only by
quadratures, since a Ricatti equation becomes linear once we know
a solution; with the {\it Bianchi Permutability Theorem} (BPT) of
1892 the B transformation can be iterated using only algebraic
computations and derivatives once we know all B transforms of the
seed) and found the spectral family of a CGC $-1$ surface (similar
to Bonnet's spectral family of CGC $1$ and {\it constant mean
curvature} (CMC) surfaces); the B transformation of CGC $-1$
surfaces is a conjugation of Bianchi's complementary
transformation with Lie's transformation.

Note that although the B transformation of the pseudo-sphere does
not differ much from Bianchi's complementary transformation from
an analytic point of view, it is essentially different from a
geometric point of view, as for the later the reflected
distributions of facets coincide (the collapsed leaves are the
rulings of the isotropic cone), while for the former the reflected
distributions of facets are different (the collapsed leaves are
the rulings of the two ruling families of a confocal
pseudo-sphere).

This collapsing ansatz allows us to simplify the denomination {\it
B transformation of surfaces applicable to quadrics} to {\it B
transformation of quadrics}; Bianchi's B transformation of
quadrics is just the metric-projective generalization of Lie's
point of view on the B transformation of the pseudo-sphere (one
replaces 'pseudo-sphere' with 'quadric' and 'circle' with
'conic').

Bianchi considered the most general form of a B transformation as
the focal surfaces (one transform of the other) of a {\it
Weingarten} (W) congruence (congruence upon whose two focal
surfaces the asymptotic directions correspond; equivalently the
second fundamental forms are proportional). Note that although the
correspondence provided by the W congruence does not give the
applicability correspondence, the B transformation is the tool
best suited to attack the deformation problem via transformation,
since it provides correspondence of the characteristics of the
deformation problem (according to Darboux these are the asymptotic
directions), it is directly linked to the infinitesimal
deformation problem (Darboux proved that infinitesimal
deformations generate W congruences and Guichard proved the
converse: there is an infinitesimal deformation of a focal surface
of a W congruence in the direction normal to the other focal
surface) and it admits a version of the BPT for its second
iteration.

First we consider general $3$-dimensional IRDF and obtain the {\it
integrability condition} (IC).

There are several basic facts of the $B_z$ transformation of
quadrics that can be considered separately or partially grouped
for $3$-dimensional IRDF: the TC requirement (tangential
distributions of facets) that facets are centered on tangent
planes, the symmetric TC requirement (facets further pass through
the origin of the tangent planes), W congruence property (which
assumes the symmetric TC), {\it applicability correspondence of
leaves of a general nature} (ACLGN) (independent of the shape of
the seed), collapsing ansatz of leaves, etc; for example Bianchi
\cite{B3} obtained the (isotropic) singular $B_z$ transformation
of quadrics with $x_z$ plane by assuming the symmetric TC with
collapsing to curves ansatz and $x_z$ plane (note that the leaves
are applicable to a quadric different from that of the
applicability of the seed for the isotropic singular $B_z$
transformation of quadrics except when $x_z$ is an isotropic plane
of a pencil of (isotropic) planes).

We derive the IC of a $3$-dimensional IRDF with the symmetric TC
and collapsing to curves ansatz (Bianchi's original assumption
from \cite{B3} except for $x_z$ being a plane); in this case the
IC is a relation which depends from $x_0$ on $x_0$ and its normal
$N_0$ and from the auxiliary surface $x_z$ on $x_z$ and its first
and second mixed derivatives in the special coordinates that give
the degenerate leaves and must be a consequence of the TC (since
the TC establishes a functional relationship between the four
parameters of $x_0$ and $x_z$, leaving only three of them
independent to generate the $3$-dimensional IRDF, any other needed
functional relationship between these four parameters must be a
consequence of the TC); a simple analysis of the prerequisite
showed that the {\it symmetric} assumption can be dropped.

In this case we found that the ACLGN does not impose any other
condition (it is a consequence of the IC) and thus all B
transformations with defining surface are pertinent to the
deformation problem (only the W congruence property required by
Bianchi's definition of a B transformation remained to be proved).

By assuming that $x_z$ is a quadric doubly ruled by leaves the
influence from $x_z$ in the IC dissolves completely with the TC
used only to symmetrize the IC; the remaining IC depends only on
$x_0$ and $N_0$ and leads to $x_0$ being a quadric confocal to
$x_z$ (again isotropic developables circumscribed to conics appear
as singular cases in our discussion).

Conversely, when the TC only symmetrizes the IC we get $x_z$
((isotropic) singular) quadric (including in a general sense
(isotropic) plane of certain pencils of (isotropic) planes) doubly
ruled by leaves, $x_0$ quadric confocal to $x_z$ and we capture
all ((isotropic) singular) finite B transformations of quadrics;
thus we complete Bianchi's approach \cite{B3} to isotropic planes
(Bianchi may have been aware of these isotropic cases, since some
of the B transformations in these cases are inverses of his B
transformations).

In trying to prove the rigidity of the B transformation of
quadrics we were led to discuss the case $x_0$ quadric, in which
case $x_z$ must be as in the case when the TC only symmetrizes the
IC.

Our IC of $3$-dimensional tangential IRDF with collapsing to
curves ansatz is not suitable to making hypothesis on the
auxiliary surface $x_z$ for  $x_z$ (isotropic) plane or quadric,
since it is very strong in assuming also the curves given by
degenerate leaves.

Bianchi's approach from \cite{B3} assumes collapsing ansatz to
$x_z$ plane only to introduce a (linear) parameter $w$ in a
particular $3$-dimensional IRDF with the symmetric TC in order to
derive its IC; he was able to prove the more general statement of
rigidity of the (isotropic) singular B transformation with $x_z$
plane by assuming that $x_z$ is a plane (with no hypothesis made
on the degenerate leaves).

Having Bianchi's approach in mind, we impose the TC on the IC of
the general $3$-dimensional IRDF and obtain some simplification of
the IC; by imposing the symmetric TC we get further
simplifications of the IC similar to those of Bianchi \cite{B3}.

In trying to see what conditions are obtained by imposing the W
congruence property (a necessary and sufficient condition to
obtain B transformation according to Bianchi's definition; the
only remaining step to obtain B transformation pertinent to the
deformation problem is to impose the ACLGN) on a $3$-dimensional
IRDF we found that there are no such conditions: the W congruence
property is a consequence of the IC of a $3$-dimensional IRDF with
the symmetric TC (in fact the W congruence property individuates
the first compatibility condition of such an IC). According to a
theorem of Ribacour $x^0$ and $x^1$ are the focal surfaces of a W
congruence if and only if
$K(x^0)K(x^1)=\frac{\sin^4(\be)}{\mathbf{d}^4}$, where $\be$
(respectively $\mathbf{d}$) is the angle (the distance) between
the tangent planes at the corresponding points (the corresponding
points); all these quantities except $K(x^1)$ are preserved for
different rolled leaves, so we get {\it same Gau\ss\ curvature
correspondence of leaves of a general nature}, which is close to
ACLGN.

Thus we obtain an important theorem:

\begin{theorem}\label{th:th1}
The necessary and sufficient condition for a $2$-dimensional
linear element to admit B\"{a}cklund transformation of its
deformations*\footnote{* Bianchi quotes in \cite{B3} two articles
of B\"{a}cklund's from 1914 and 1916 on Bianchi's B transformation
of quadrics, so B\"{a}cklund himself may have some other results
on the deformation problem via the B transformation, but
unfortunately we are not familiar with German.} is
(\ref{eq:dissymTC1}).
\end{theorem}

By imposing post-priori the collapsing to curves ansatz on the IC
of a $3$-dimensional IRDF with the symmetric TC we were able to
show that $x_0$ must be a quadric confocal to $x_z$ when $x_z$ is
a quadric; thus by the previous argument $x_z$ must be doubly
ruled by degenerate leaves (the case $x_z$ plane is due to Bianchi
\cite{B3}; the same argument applies to the isotropic case and
from Bianchi's argument we only need the part up to obtaining
$x_z$ doubly ruled by degenerate leaves).

Using both approaches we were finally able to prove:
\begin{theorem}\label{th:th2}
The only B\"{a}cklund transformation with defining surface is
Bianchi's B\"{a}cklund transformation of quadrics.
\end{theorem}

It remains to consider in detail the B transformation of
deformations of abstract linear elements (without defining
surface)and with ACLGN.

By requiring that the leaves are applicable to certain surfaces
with ACLGN (independent of the shape of $x$) and since the
distribution of facets into leaves changes with the shape of $x$
we get the fact that all leaves of all rolled distributions of
facets must be applicable to the same surface $y$; further we get
the TC and the necessary algebraic conditions for such a
configuration to exist (it depends on three arbitrary constants)
are satisfied for the symmetric TC.

By imposing compatibility conditions on these algebraic conditions
we must get the space of solutions depending on two functions of a
variable for the deformations of $y$; again these conditions may
be satisfied by $3$-dimensional IRDF with the symmetric TC or
further compatibility conditions may be imposed on the linear
element of $x$ which may make the family of such surfaces
depending on constants or on a function of one variable and
constants.

The change from a facet $\mathcal{F}$ to an infinitesimally close
facet $\mathcal{F}'$ via infinitesimal rolling in an arbitrary
common tangential infinitesimal direction $\del$ can be reversed
by considering $-\del$ acting on $\mathcal{F}'$; thus this
infinitesimal law corresponds via discretization to the
prescription of the finite ACLGN (correspondence between the leaf
$x^1$ and a surface $x_0^1\subset x_0$ or between $x^1$ and a
surface $x_0^1\subset y_0$ with $y_0$ another (defining) surface;
according to Bianchi \cite{B3} there are B transformations such
that the seed and the leaf are applicable to different quadrics)
and the inversion of the B transformation (the seed and the leaf
exchange places; note that in this case the leaf applicable to the
(defining) surface $x_0$ is required to belong to an
$1$-dimensional family, which is certainly the case if the
original leaf can play the r\^{o}le of seed (that is it admits B
transformation)).

In dealing with iterations of infinitesimal rollings in arbitrary
tangential infinitesimal directions we are led to consider their
(partial) discretizations; these are realized by M\"{o}bius
configurations.

A tetrahedron consists of $2^2$ points and $2^2$ planes, each
point (plane) belonging to (containing) $2+1$ planes (points).
M\"{o}bius considered configurations $\mathcal{M}_3$ of two
tetrahedra inscribed one into the other, that is configurations of
$2^3$ points and $2^3$ planes, each point (plane) belonging to
(containing) $3+1$ planes (points); two $\mathcal{M}_3$
configurations inscribed one into the other give raise to a
configuration of $2^4$ points and $2^4$ planes, each point (plane)
belonging to (containing) $4+1$ planes (points). Therefore Bianchi
(\cite{B5},Vol {\bf 5},(117)) calls a configuration of $2^n$
points and $2^n$ planes such that any point (plane) belongs to
(contains) $n+1$ planes (points) a {\it M\"{o}bius configuration}
$\mathcal{M}_n$.

For the discussion up to the {\it third iteration of the tangency
configuration} (TITC) we consider only the (singular) B
transformation of quadrics, when the seed and the leaf are
applicable via the Ivory affinity to the same quadric.

If we compose the inverse of the rolling of $x_0^1\subset x_0$ on
the leaf $x^1$ with the rolling of $x_0^0$ on the seed $x^0$, then
we get a rigid motion at the level of the static picture with the
defining surface $x_0$ and auxiliary surface $x_z$ (the {\it rigid
motion provided by the Ivory affinity} (RMPIA); this led Bianchi
to discover the {\it applicability correspondence provided by the
Ivory affinity} (ACPIA) in \cite{B1}); thus the prescription of
the applicability correspondence becomes of a general nature
(independent of the shape of the seed $x^0$) and becomes, when the
leaves collapse, a correspondence between $x_z$ and $x_0$. In this
correspondence the special coordinates on $x_z$ (rulings given by
the collapsed leaves) correspond to special coordinates on $x_0$,
thus providing special coordinates for the whole problem.

We have the symmetry of the TC: four $1$-M\"{o}bius configurations
$\mathcal{M}_1$, one taken to two others via a reflection in a
facet or a RMPIA; the symmetry of the TC implies the inversion of
the B transformation and existence of RMPIA in the TC (note that
a-priori the RMPIA does not require the TC) implies prescription
of the applicability correspondence as a law of a general nature
(independent of the shape of the seed).

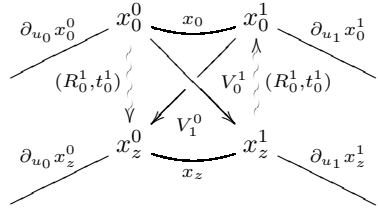
\begin{center}
$\xymatrix@!0{&&x_0^0\ar@{-}[drdr]\ar@/_/@{-}[rr]^{x_0}\ar@{~>}[dd]_{(R_0^1,t_0^1)}&&
x_0^1\ar@{<~}[dd]^{(R_0^1,t_0^1)}&&\\
\ar@{-}[urr]^{\pa_{u_0}x_0^0}&&&\ar[dl]^>>>>{V_1^0}\ar[dr]^>>>>>{V_0^1}&&&
\ar@{-}[ull]_{\pa_{u_1}x_0^1}\\&&
x_z^0\ar@{-}'[ur][urur]\ar@/_/@{-}[rr]_{x_z}&&x_z^1&\\\ar@{-}[urr]^{\pa_{u_0}x_z^0}&&&
&&&\ar@{-}[ull]_{\pa_{u_1}x_z^1}}$
\end{center}
$$(R_j,t_j)(x_0^j,dx_0^j):=(R_jx_0^j+t_j,R_jdx_0^j)=(x^j,dx^j),\
(R_j,t_j)\subset\mathbf{O}_{3}(\mathbb{C})\ltimes\mathbb{C}^{3},\
j=0,1,$$ $$(R_0^1,t_0^1)=(R_1,t_1)^{-1}(R_0,t_0),\
V_0^1:=x_z^1-x_0^0,\ V_1^0:=x_z^0-x_0^1,\
(V_0^1)^TN_0^0=0\Leftrightarrow(V_1^0)^TN_0^1=0,$$ $R_0^1[V_0^1\ \
\pa_{u_1}x_z^1\ \ \pa_{u_0}x_0^0]=[-V_1^0\ \ \pa_{u_1}x_0^1\ \
\pa_{u_0}x_z^0],\ (V_0^1)^TN_0^0=0\Rightarrow
R_0^1(I_3-2N_0^0(N_0^0)^T)\pa_{v_1}x_z^1=\pa_{v_1}x_0^1$.

For these arguments we need only a spectral parameter $z$ (a
$3$-dimensional IRDF); for the remaining arguments we need $z$ to
vary.

The issue of discretizing the commutation of the composition of
infinitesimal rolling in an arbitrary tangential direction $\del$,
followed by infinitesimal rolling in an arbitrary tangential
direction $\del'$: $\del\del'=\del'\del$ (equivalent to the
symmetry of the difference of the second fundamental forms of the
rolling surfaces) in one of $\del,\del'$ (in which case the other
remains infinitesimal) leads to finding the differential system
subjacent to the B transformation (Ricatti equation) and proving
the applicability correspondence; in both of $\del,\del'$ leads to
algebraic computations of the {\it second iterated tangency
configuration} (SITC) (four $2$-M\"{o}bius configurations
$\mathcal{M}_2$, each taken to two others via a RMPIA in the TC
and to the fourth via a RMPIA without the TC): for $x_{z_1}^0,\
x_{z_2}^3\in T_{x_0^1}x_0$ in the pencil of planes containing
$x_{z_2}^0,\ x_{z_1}^3$ there are two planes tangent to $x_0$ at
$x_0^2$; thus we have two choices of $x_0^2$ according to the
rulings of $x_{z_1}^0,\ x_{z_2}^3$ belonging to the same or
different ruling families.

Once the SITC is established to be valid, one can let one of
$\del,\del'$ be infinitesimal and obtain the differential system
subjacent to the B transformation and proof of the applicability
correspondence, so the SITC is sufficient to imply the
differential system subjacent the B transformation and the proof
of the applicability correspondence (that is the first {\it moving
M\"{o}bius configuration}). Note that the SITC is equivalent to
the rulings at $x_{z_1}^0,\ x_{z_2}^3$ cutting the segment
$[x_0^1\ (R_3^0,t_3^0)x_0^2]$ with cross-ratio $\frac{z_1}{z_2}$
and it boils down to a homography between four rulings and having
the symmetries of the square (the Bianchi quadrilateral); the
rulings of opposite vertices belong simultaneously to the same or
different ruling families.

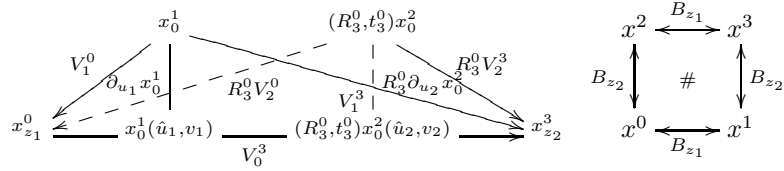
\begin{center}
$\xymatrix{&
^{x_0^1}\ar[dl]_{V_1^0}\ar@{-}[d]_{\pa_{u_1}x_0^1}\ar[drr]_<<<<<<<<<<<<<<<<<<<<<<<{V_1^3}&
^{(R_3^0,t_3^0)x_0^2}\ar@{-->}[dll]^<<<<<<<<<<<<{R_3^0V_2^0}
\ar@{--}[d]^{R_3^0\pa_{u_2}x_0^2}
\ar[dr]^{R_3^0V_2^3}&\\^{x_{z_1}^0}\ar@{-}[r]&^{x_0^1(\hat
u_1,v_1)}\ar@{-}[r]_<<<<<{V_0^3}&^{(R_3^0,t_3^0)x_0^2(\hat
u_2,v_2)}\ar[r]&^{x_{z_2}^3}}$ $\xymatrix{\ar@{}[dr]|{\#}x^2\ar@{<->}[d]_{B_{z_2}}\ar@{<->}[r]^{B_{z_1}}&x^3\ar@{<->}[d]^{B_{z_2}}\\
x^0\ar@{<->}[r]_{B_{z_1}}&x^1}$
\end{center}
$$(R_j,t_j)(x_0^j,dx_0^j)=(x^j,dx^j),\ j=0,1,2,3,$$
$$(R_j^k,t_j^k)=(R_k,t_k)^{-1}(R_j,t_j),\
(j,k)=(0,1),(0,2),(1,3),(2,3),$$
$$(R_0^1,t_0^1)(R_2^0,t_2^0)=(R_3^0,t_3^0)=(R_3^1,t_3^1)(R_2^3,t_2^3),\
(R_1^0,t_1^0)(R_3^1,t_3^1)=(R_2^1,t_2^1)=(R_2^0,t_2^0)(R_3^2,t_3^2).$$
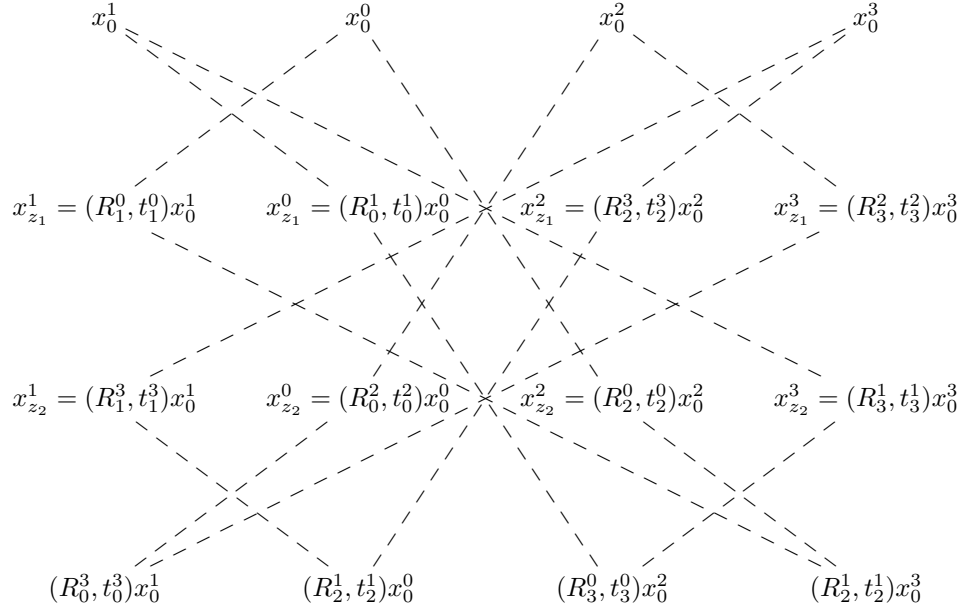
\begin{center}
$\xymatrix@!0{x_0^1\ar@{--}[dddrrrr]\ar@{--}[ddddddrrrrrrrrrrrr]&&&&
x_0^0\ar@{--}[dddllll]\ar@{--}[ddddddrrrr]&&&&
x_0^2\ar@{--}[ddddddllll]\ar@{--}[dddrrrr]&&&&
x_0^3\ar@{--}[ddddddllllllllllll]\ar@{--}[dddllll]\\
\\  \\
x_{z_1}^1=(R_1^0,t_1^0)x_0^1\ar@{--}[ddddddrrrrrrrrrrrr]&&&&
x_{z_1}^0=(R_0^1,t_0^1)x_0^0\ar@{--}[ddddddrrrr]&&&&
x_{z_1}^2=(R_2^3,t_2^3)x_0^2\ar@{--}[ddddddllll]&&&&
x_{z_1}^3=(R_3^2,t_3^2)x_0^3\ar@{--}[ddddddllllllllllll]\\
\\  \\
x_{z_2}^1=(R_1^3,t_1^3)x_0^1\ar@{--}[dddrrrr]&&&&
x_{z_2}^0=(R_0^2,t_0^2)x_0^0\ar@{--}[dddllll]&&&&
x_{z_2}^2=(R_2^0,t_2^0)x_0^2\ar@{--}[dddrrrr]&&&&
x_{z_2}^3=(R_3^1,t_3^1)x_0^3\ar@{--}[dddllll]\\
\\   \\
(R_0^3,t_0^3)x_0^1&&&&(R_2^1,t_2^1)x_0^0&&&&(R_3^0,t_3^0)x_0^2&&&&(R_2^1,t_2^1)x_0^3}$
\end{center}

The issue of discretizing the commutation of the composition of
infinitesimal rolling in an arbitrary tangential direction $\del$,
followed by infinitesimal rolling in an arbitrary tangential
direction $\del'$, followed by infinitesimal rolling in an
arbitrary tangential direction $\del''$: $\del\del'\del''
=\del'\del''\del$ (equivalent to the difference of the G-CMP
equations of the rolling surfaces or the flatness of the flat
connection form) in one of $\del,\del',\del''$ (in which case the
other two remain infinitesimal) leads to the existence of the B
transformation (equivalently the complete integrability of the
differential system subjacent to the B transformation or the
integrability of the considered $3$-dimensional rolled
distributions of facets); in two of $\del,\del',\del''$ (in which
case the other remains infinitesimal) leads to the BPT and in all
three $\del,\ \del',\ \del''$ leads to the algebraic computations
of the TITC (eight $3$-M\"{o}bius configurations $\mathcal{M}_3$,
each taken to three others via a RMPIA in the TC and to three
others via a RMPIA without the TC). Once the TITC is established
to be valid, one can let one of $\del,\del',\del''$ be
infinitesimal and obtain the BPT (that is the second {\it moving
M\"{o}bius configuration}) or one can let two of
$\del,\del',\del''$ be infinitesimal and obtain the full theory
(existence of B transformation). Note that the TITC uses the
cross-ratio properties with
$\frac{z_1}{z_2},\frac{z_2}{z_3},\frac{z_3}{z_1}$ and the Menelaos
theorem; it is equivalent to a homography between four rulings of
corresponding sameness ruling families and having the symmetries
of the cube; the four rulings involved are from the vertices of a
regular tetrahedron and thus are in symmetric relationship one to
the other if they belong to the same ruling family.

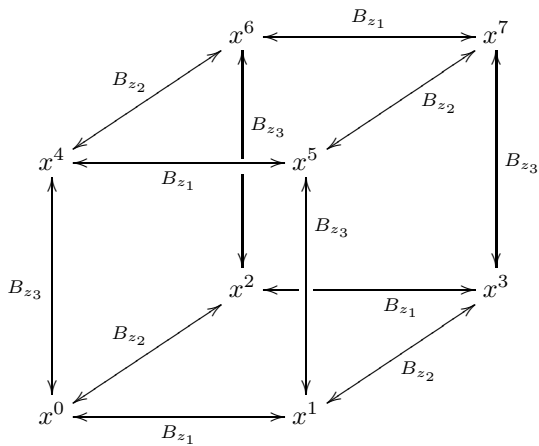
\begin{center}
$\xymatrix@!0{&&&x^6\ar@{<->}[rrrr]^{B_{z_1}}&&&&x^7\\
&&&&&&&\\
x^4\ar@{<->}[uurrr]^{B_{z_2}} \ar@{<->}[rrrr]_{B_{z_1}}&&&&
x^5\ar@{<->}[uurrr]_>>>>>>>>>{B_{z_2}}&&&\\
&&&&&&&\\
&&&x^2\ar@{<->}'[r][rrrr]_{B_{z_1}}
\ar@{<->}'[uu][uuuu]_<<<<<{B_{z_3}}&&&&
x^3\ar@{<->}[uuuu]_{B_{z_3}}\\
&&&&&&&\\
x^0\ar@{<->}[rrrr]_{B_{z_1}}\ar@{<->}[uurrr]^{B_{z_2}}
\ar@{<->}[uuuu]^{B_{z_3}}&&&& x^1\ar@{<->}[uurrr]_{B_{z_2}}
\ar@{<->}[uuuu]_<<<<<<<<<<<<<<<<<<<<<<{B_{z_3}}&&&}$
\end{center}

Existence of $3$-M\"{o}bius configurations $\mathcal{M}_3$ as in
the TITC implies existence of arbitrary $n$-M\"{o}bius
configurations $\mathcal{M}_n$ for arbitrary $n$ iterations of the
TC, as expected (the G-CMP equations involve $3$ derivatives and
are necessary and sufficient for the characterization of a
surface); thus the algebraic computations of the TITC are
necessary and sufficient for a complete description of the B
transformation (some interesting exterior algebra type formulae
but no relevant information can be found beyond the TITC).

Note also that the SITC and the TITC behave well with respect to
totally real considerations, which implies the B transformation of
totally real quadrics.

These considerations can be extended to include the isotropic
singular B transformations (the RMPIA is replaced with a suitable
rigid motion; one cannot construct the isotropic singular B
transformation using the Ivory affinity since one cannot take
square roots of symmetric matrices with isotropic kernels and for
quadrics of revolution the Ivory affinity has image only the
singular part of the singular isotropic confocal quadric (the axis
of revolution)); one gets the SITC and the TITC (and thus the BPT
and higher moving M\"{o}bius configurations) with certain
iterations exchanging the quadric of applicability (most isotropic
singular B transformations have seed and leaf applicable to
different quadrics).

There is another application of isotropic developables
circumscribed to conics, which leads to the {\it Hazzidakis} (H)
transformation of quadrics (an involutory algebraic transformation
of deformations of quadrics which commutes with the B
transformation): since homographies preserve tangency and its
order, they preserve asymptotes (whose osculating planes are the
tangent planes of the surface), so it takes developables (on which
the two asymptotic directions coincide) to developables; the
homography $H$ takes the isotropic developable circumscribed to
$H^{-1}(C(\infty))$ to the isotropic developable circumscribed to
$H(C(\infty))$; for most homographies $H$ the isotropic
developable circumscribed to $H^{-1}(C(\infty))$ generates a
family of confocal quadrics.

Now by the Chasles-Jacobi result (Jacobi proved that the tangents
of a geodesic on a quadric remain tangent to a confocal quadric
and Chasles the converse: the congruence of common tangents to two
confocal quadrics is normal and envelopes geodesics on the two
quadrics) $H$ takes geodesics on quadrics to geodesics on
quadrics, so it provides pairs of quadrics {\it conjugate in
deformation} (two non-flat non-homothetic surfaces with a
correspondence of asymptotes and of all virtual asymptotes
(coordinates that become asymptotes on a deformation); it is
enough to have correspondence of asymptotes and geodesics); only
quadrics can be conjugate in deformation (Bianchi proved the
statement for surfaces of revolution and Servant the general
statement).

By use of the H transformation singular finite B transformations
are taken to singular finite B transformations (Bianchi) or
Calapso's singular $B_{\infty}$ transformation of {\it quadrics
with center} (QC) (see \cite{C}; note that in this case one knows
the $B_\infty$ transform only infinitesimally (the first and
second fundamental forms of the leaf) and we don't have the W
congruence); this behavior may be generalized to isotropic
singular B transformations taken to isotropic singular B
transformations or $B_{\infty}$ transformations of {\it
(isotropic) quadrics without center} (I)QWC) (there is some
partial evidence in favor of existence of such $B_{\infty}$
transformations).

\section{Rolling surfaces}

The study of the rolling problem was initiated by Ribacour and has
been extensively pursued in Bianchi \cite{B4},(\cite{B5},Vol {\bf
7}) and Darboux \cite{D1}.

Let $(u,v)\ni D$ with $D$ a domain in $\mathbb{R}^2,\
\mathbb{C}\times\mathbb{R}$ or $\mathbb{C}^2$ and $x:D\mapsto
\mathbb{C}^3$ be a surface.

For $\om_1,\ \om_2\ \mathbb{C}^3$-valued $1$-forms on $D$ and
$a,b\in\mathbb{C}^3$ we have
\begin{eqnarray}\label{eq:che}
a^T\om_1\wedge b^T\om_2=((a\times
b)\times\om_1+b^T\om_1a)^T\wedge\om_2=(a\times
b)^T(\om_1\times\wedge\om_2)+b^T\om_1\wedge a^T\om_2;\nonumber\\
\mathrm{in\ particular}\ a^T\om\wedge b^T\om=\frac{1}{2}(a\times
b)^T(\om\times\wedge\om).
\end{eqnarray}
Since both $\times$ and $\wedge$ are skew-symmetric, we have
$2\om_1\times\wedge\om_2=\om_1\times\om_2+\om_2\times\om_1=2\om_2\times\wedge\om_1$.

Consider the scalar product $<,>$ on $\mathbf{M}_3(\mathbb{C})$:
$<X,Y>:=\frac{1}{2}$tr$(X^TY)$. We have the isometry

$\al :\mathbb{C}^3\rightarrow\mathbf{o}_3(\mathbb{C}),\ \al(\begin{bmatrix}x^1\\
x^2\\
x^3
\end{bmatrix})=\begin{bmatrix}
0&-x^3&x^2\\
x^3&0&-x^1\\
-x^2&x^1&0
\end{bmatrix},\
x^Ty=<\al(x),\al(y)>=\frac{1}{2}$tr$(\al(x)^T\al(y)),\\
\al(x\times y)=[\al(x),\al(y)]=\al(\al(x)y)=yx^T-xy^T,\
\al(Rx)=R\al(x)R^{-1},\ x,y\in\mathbb{C}^3,\
R\in\mathbf{O}_3(\mathbb{C})$.

Let $x\subset\mathbb{C}^3$ be a surface non-rigidly applicable to
a surface $x_0\subset\mathbb{C}^3$:
\begin{eqnarray}\label{eq:roll}
(x,dx)=(R,t)(x_0,dx_0):=(Rx_0+t,Rdx_0),
\end{eqnarray}
where $(R,t)$ is sub-manifold in
$\mathbf{O}_3(\mathbb{C})\ltimes\mathbb{C}^3$ (in general surface,
but it is a curve if $x_0,\ x$ are ruled and the rulings
correspond under the applicability or a point if $x_0,x$ differ by
a rigid motion). The sub-manifold $R$ gives the rolling of $x_0$
on $x$, that is if we rigidly roll $x_0$ on $x$ such that points
corresponding under the applicability will have the same
differentials, $R$ will dictate the rotation of $x_0$; the
translation $t$ will satisfy $dt=-dRx_0$.

For $(u,v)$ parametrization on $x_0,x$ and outside the locus of
isotropic (degenerate) linear element of $x_0,x$ we have
$N_0:=\frac{\pa_ux_0\times\pa_vx_0}{|\pa_ux_0\times\pa_vx_0|},\
N:=\frac{\pa_ux\times\pa_vx}{|\pa_ux\times\pa_vx|}$ respectively
positively oriented unit normal fields of $x_0,x$ and $R$ is
determined by $R=[\pa_ux\ \ \pa_vx\ \ N][\pa_ux_0\ \ \pa_vx_0\ \
\det(R)N_0]^{-1}$; we take $R$ with $\det(R)=1$; thus the rotation
of the rolling with the other face of $x_0$ (or on the other face
of $x$) is $R':=R(I-2N_0N_0^T)=(I-2NN^T)R,\ \det(R')=-1$.

Therefore $\mathbf{O}_3(\mathbb{C})\ltimes\mathbb{C}^3$ acts on
$2$-dimensional integrable distributions of facets $(x_0,dx_0)$ in
$T^*(\mathbb{C}^3)$ as: $(R,t)(x_0,dx_0)=(Rx_0+t,Rdx_0)$; a
rolling is a sub-manifold
$(R,t)\subset\mathbf{O}_3(\mathbb{C})\ltimes\mathbb{C}^3$ such
that $(R,t)(x_0,dx_0)$ is still integrable.

We have:
\begin{eqnarray}\label{eq:secoi}
R^{-1}dRN_0=R^{-1}dN-dN_0.
\end{eqnarray}
In order to preserve the classical notation $d^2$ for the
tensorial (symmetric) second derivative we shall use $d\wedge$ for
the exterior (antisymmetric) derivative. Applying the
compatibility condition $d\wedge$ to (\ref{eq:roll}) we get:
\begin{eqnarray}\label{eq:comp}
R^{-1}dR\wedge dx_0=0,\ dRR^{-1}\wedge dx=0.
\end{eqnarray}
Applying $R^{-1}d$ to (\ref{eq:roll}) we get
\begin{eqnarray}\label{eq:roll1}
R^{-1}d^2x=R^{-1}dRdx_0+d^2x_0.
\end{eqnarray}
Since $R^{-1}dR$ is skew-symmetric and using (\ref{eq:comp}) we
have
\begin{eqnarray}\label{eq:dx0}
dx_0^TR^{-1}dRdx_0=0.
\end{eqnarray}
From (\ref{eq:dx0}) for $a\in\mathbb{C}^3$ we get
$R^{-1}dRa=R^{-1}dR(a^{\bot}+a^{\top})=a^TN_0R^{-1}dRN_0-a^TR^{-1}dRN_0N_0=\om\times
a,\ \om:=N_0\times R^{-1}dRN_0=^{(\ref{eq:secoi})}(\det
R)R^{-1}(N\times dN)-N_0\times dN_0=R^{-1}(N\times dN)-N_0\times
dN_0$. Thus $R^{-1}dR=\al(\om)$ and $\om$ is flat connection form
in $T^*x_0$:
\begin{eqnarray}\label{eq:om}
d\wedge\om+\frac{1}{2}\om\times\wedge\om=0,\ \om\times\wedge
dx_0=0,\ (\om)^\perp=0.
\end{eqnarray}
With $s:=N_0^T(R^{-1}d^2x-d^2x_0)
=^{(\ref{eq:roll1})}N_0^T(\om\times
dx_0)=s_{11}du^2+s_{12}dudv+s_{21}dvdu+s_{22}dv^2$ the difference
of the second fundamental forms of $x,\ x_0$ we have
\begin{eqnarray}\label{eq:omjk}
\om=\frac{s_{12}\pa_ux_0-s_{11}\pa_vx_0}{|\pa_ux_0\times\pa_vx_0|}du+
\frac{s_{22}\pa_ux_0-s_{21}\pa_vx_0}{|\pa_ux_0\times\pa_vx_0|}dv;
\end{eqnarray}
($\om\times\wedge dx_0=0$ is equivalent to $s_{12}=s_{21}$;
$(d\wedge\om)^\perp+\frac{1}{2}\om\times\wedge\om=0,\
(d\wedge\om)^\top=0$ respectively encode the difference of the
G-CMP equations of $x_0$ and $x$).

Using $\frac{1}{2}dN_0\times\wedge
dN_0=K|\pa_ux_0\times\pa_vx_0|N_0du\wedge dv,\ K$ being the
Gau\ss\ curvature we get $dN_0\times\wedge
dN_0=R^{-1}(dN\times\wedge dN)=^{(\ref{eq:secoi})}(\om\times
N_0+dN_0)\times\wedge(\om\times N_0+dN_0)=dN_0\times\wedge
dN_0+2(\om\times N_0)\times\wedge dN_0+\om\times\wedge\om$; thus
\begin{eqnarray}\label{eq:omom}
\frac{1}{2}\om\times\wedge\om=dN_0^T\wedge\om N_0.
\end{eqnarray}
Note also
\begin{eqnarray}\label{eq:om'}
\om'=N_0\times {R'}^{-1}dR'N_0=-\om-2N_0\times dN_0
\end{eqnarray}
and
\begin{eqnarray}\label{eq:aom}
\ \ \ \ a^T\wedge\om=0,\ \forall\om\ \mathrm{satisfying\
(\ref{eq:om})\ for}\ a\ 1-\mathrm{form}\Rightarrow\ a^T\odot
dx_0:=\frac{a^Tdx_0+dx_0^Ta}{2}=0.
\end{eqnarray}
Note that the converse $a^T\odot dx_0^0=0,\ a\ 1-$form
$\Rightarrow a^T\wedge\om=0,\forall\om$ satisfying (\ref{eq:om})
is also true.

\subsection{Isotropic developables}

For any two curves $c_1(v),c_2(v)$ the developable circumscribed
to them is $(u,v)\mapsto c_1(v)+u[c_2(f(v))-c_1(v)]$, where $f(v)$
is determined from $[c_2(f(v))-c_1(v)]^T[c_1'(v)\times
c_2'(f(v))]=0$, that is $c_2(f(v))-c_1(v)$ belongs to tangent
planes of $c_1$ and $c_2$; isotropic developables are the
developables circumscribed to a finite curve and $C(\infty)$.

With $e_1,e_2,e_3,\ e_j^Te_k=\del_{jk}$ the standard basis of
$\mathbb{C}^3$ and $f_1:=\frac{e_1-ie_2}{\sqrt{2}}$ the standard
isotropic vector we have $Y(v):=-v^2f_1+2\bar f_1+2ve_3$ the
standard parametrization of the rulings of the isotropic cone and
the isotropic developable circumscribed to $c(v)$ is $(u,v)\mapsto
uY(f(v))+c(v),\ c'(v)^TY(f(v))=0$ (there are two choices of $f(v)$
except for $c'(v)$ isotropic).

If $x_0\subset\mathbb{C}^3$ is a surface with degenerate linear
element, then its tangent planes are isotropic and we can take the
curves whose tangents are isotropic as the curves $v=$ct; thus
$\pa_ux_0=a(u,v)Y(w(u,v)),\
\pa_vx_0=b(u,v)Y(w(u,v))+c(u,v)Y'(w(u,v))$; after a change of
coordinates we get $\pa_ux_0=a(u,v)Y(v)$ (so after another change
of coordinates $x_0(u,v)=uY(v)+c(v),\ Y(v)^Tc'(v)=0$) or
$\pa_ux_0=a(u,v)Y(u),\ \pa_vx_0=b(u,v)Y(u)+c(u,v)Y'(u)$ (in this
case from $\pa_{uv}^2x_0=\pa_vaY(u)$ we get $b=c=0$, a
contradiction) or $x_0=uY(w)+vY'(w),\ w=$ct.

\section{$3$-dimensional integrable rolling distributions of facets}

Assume that we have a $3$-dimensional distribution of facets
$$(p,P)=(p(u,v,w),P(u,v,w)),\ p\in P,\ du\wedge dv\wedge dw\neq 0$$
in $\mathbb{C}^3$ with normal fields
$m=m(u,v,w)\subset\mathbb{C}^3\setminus\{0\},\ m\bot P$.

With $\ti d\cdot:=\pa_u\cdot du+\pa_v\cdot dv+\pa_w\cdot
dw=d\cdot+\pa_w\cdot dw$ if the distribution of facets is
integrable, then along the leaves we have
\begin{eqnarray}\label{eq:inte}
0=m^T\ti dp=m^T(\pa_updu+\pa_vpdv+\pa_wpdw)=m^T(dp+\pa_wpdw).
\end{eqnarray}
Assuming $m^T\pa_wp\neq 0$, applying the compatibility condition
$m^T\pa_wp\ti d\wedge$ to (\ref{eq:inte}) and using the equation
itself we get the IC $m^T\pa_wp\neq 0,\
(\pa_wm^Tdp-dm^T\pa_wp)\wedge m^Tdp+m^T\pa_wpdm^T\wedge dp=0$, or:
\begin{eqnarray}\label{eq:inteco}
\ \ \ \ \ \ \ m^T\pa_wp\neq 0,\ (dp\times\pa_wp)^T\wedge(m\times
dm)+\frac{1}{2}(\pa_wm\times m)^T(dp\times\wedge dp)=0,\ du\wedge
dv\wedge dw\neq 0
\end{eqnarray}
(in order to get the $1$-dimensional family of leaves $c=$ct from
the integration of (\ref{eq:inte}) ($w=w(u,v,c)$ for
$m^T\pa_upm^T\pa_vp\neq 0$ or $w=w(u,c)$ for $m^T\pa_up\neq0,\
m^T\pa_vp=0$ or $w=w(v,c)$ for $m^T\pa_up=0,\ m^T\pa_vp\neq 0$ or
$w=c$ for $m^T\pa_up=m^T\pa_vp=0$), (\ref{eq:inteco}) must be
identically satisfied (without imposing a functional relationship
between $u,v$ and $w$; note that the scaling of $m$ is irrelevant,
but $m$ may be isotropic, in which case the leaves are isotropic
developables)).

Since along the leaves we have $\frac{1}{2}\ti dp\times\wedge\ti
dp=[(I_3-\frac{\pa_wpm^T}{m^T\pa_wp})\pa_up]\times
[(I_3-\frac{\pa_wpm^T}{m^T\pa_wp})\pa_vp]du\wedge
dv=[(I_3-\frac{\pa_wpm^T}{m^T\pa_wp})^*]^T(\pa_up\times\pa_vp)du\wedge
dv$ and $\ker[(I_3-\frac{\pa_wpm^T}{m^T\pa_wp})^*]^T=
\mathrm{Im}(I_3-\frac{\pa_wpm^T}{m^T\pa_wp})^T=(\pa_wp)^\bot$, the
leaves are $2$-dimensional unless $\pa_wp^T(\pa_up\times\pa_vp)=0$
(note that along the leaves we need $du\wedge dv\neq 0$ in order
to preserve the $3$-dimensionality of the distribution of facets).

By symmetry in the variables $(u,v,w)$ the only remaining singular
case to discuss is when along the leaves we have
$m^Tdp=m^T\pa_wp=0,\ du\wedge dv\neq 0$, in which case
$\pa_wp^T(dp\times\wedge dp)=0$, so the centers of facets are
situated on a surface, curve or point. Using
$d(m^T\pa_wp)=\pa_w(m^Tdp)=0$ and applying the compatibility
condition $\ti d\wedge$ to $0=m^Tdp=m^T\ti dp$ we get the IC
$0=dm^T\wedge dp+(dm^T\pa_wp-\pa_wm^Tdp)\wedge dw=dm^T\wedge dp$,
or
\begin{eqnarray}\label{eq:inteco1}
m^T\pa_wp=0,\ m^Tdp=0,\ dm^T\wedge dp=0,\ du\wedge dv\wedge dw\neq
0.
\end{eqnarray}
In this case the $3$-dimensional integrable distribution of facets
is just a $2$-dimensional integrable distribution of facets (the
tangent planes of a surface, curve or point), each facet being
counted with the simple $\infty$ multiplicity of $w$ (the
dependence on $w$ is irrelevant to our problem).

We can distribute the distribution of facets along the surface
$x_0=x_0(u,v)$ (the parameters $(u,v)$ and $w$ are individuated
such that to each point of $x_0$ corresponds an $1$-dimensional
family of facets of the distribution depending on $w$; thus one
must discuss singular cases according only to the symmetry
$u\leftrightarrow v$); a-priori the distribution of facets has no
other relation to $x_0$.

By referring $V:=p-x_0,m$ to $dx_0$ and $N_0$ (\ref{eq:inteco})
becomes an equation involving the geometry of $x_0$ (depending on
the linear element and linearly on the second fundamental form; by
an application of the Gau\ss\ theorem the terms of the second
fundamental form appearing quadratically group together to give
dependence on the linear element).

The natural question thus appears wether the IC (\ref{eq:inteco})
depends only on the linear element of $x_0$ (that is we require
the cancellation of the coefficients of the (linearly appearing)
second fundamental form); equivalently if we roll $x_0$ on an
applicable surface $(x,dx)=(Rx_0+t,Rdx_0)$, then (\ref{eq:inteco})
is still satisfied if we replace $V,m,x_0$ with $RV,Rm,x$ (note
that by referring $V,m$ to $dx_0$ and $N_0$, their coefficients
may depend on $dx_0$ and $N_0$, so the derivatives of these
coefficients may depend on the second fundamental form of $x_0$;
we ignore this dependence in our considerations since the
coefficients themselves are preserved by rolling).

Note $N_0^TVdN_0+d[(I_3-N_0N_0^T)V]^TN_0N_0=(N_0\times dN_0)\times
V$ is the part of $dV$ depending (linearly) on the second
fundamental form of $x_0$, so in (\ref{eq:inteco}) we need to
consider the condition that the terms that do not depend linearly
on $N_0\times dN_0$ cancel and individuate those that depend
linearly on $N_0\times dN_0$ to be cancelled separately by
replacing $(...)^T\wedge(N_0\times dN_0)=0$ with $(...)^T\odot
dx_0=0$ (thus the use of $\om=R^{-1}(N\times dN)-N_0\times dN_0$
and (\ref{eq:aom}) becomes clear).

Equation (\ref{eq:inteco}) becomes $[([d(V+x_0)-(N_0\times
dN_0)\times V]+(N_0\times dN_0)\times
V)\times\pa_wV]^T\wedge[m\times([dm-(N_0\times dN_0)\times m]
+(N_0\times dN_0)\times m)]+\frac{1}{2}(\pa_wm\times m)^T
[([d(V+x_0)-(N_0\times dN_0)\times V]+(N_0\times dN_0)\times
V)\times\wedge([d(V+x_0)-(N_0\times dN_0)\times V]+(N_0\times
dN_0)\times V)]=0$ (note that in the reflected distribution of
facets $V':=(I_3-2N_0N_0^T)V,m':=(I_3-2N_0N_0^T)m$ the terms
depending linearly on the second fundamental form have changed
signs); by separating the terms as explained we get the IC
\begin{eqnarray}\label{eq:dis}
m^T\pa_wV\neq 0,\ 2[(d(V+x_0)-(N_0\times dN_0)\times
V)\times\pa_wV]^T\wedge[m\times(dm-(N_0\times dN_0)\times
m)]+\nonumber\\ (\pa_wm\times m)^T[(d(V+x_0)-(N_0\times
dN_0)\times V)\times\wedge(d(V+x_0)-(N_0\times dN_0)\times V)]\nonumber\\
+[\pa_w(m\times V)\times(m\times V)+m^T\pa_wVm\times
V]^TN_0N_0^T(dN_0 \times\wedge dN_0)=0,\ ([d(V+x_0)-\nonumber\\
(N_0\times dN_0)\times V]^T[m\pa_w(m\times V)-\pa_wm(m\times
V)]+\pa_wV^T(dm-(N_0\times dN_0)\times m)m\times
V\nonumber\\-m^T\pa_wV[d(m\times V)-(N_0\times dN_0)\times(m\times
V)])^T\odot dx_0=0,\ du\wedge dv\wedge dw\neq 0.
\end{eqnarray}
Thus we have a $3$-dimensional IRDF provided (\ref{eq:dis}) is
identically satisfied (without imposing a functional relationship
between $u,v$ and $w$).

In the case of tangential distributions of facets ($V^TN_0=0$)
(\ref{eq:inte}) for the rolled distribution of facets becomes
$m^T[-V^T(\om\times N_0)N_0+d(V+x_0)+\pa_wVdw]=0$; imposing the
compatibility condition $m^T\pa_wV\ti d\wedge$, using
(\ref{eq:om}) and the equation itself we get
$0=m^T\pa_wV[-d(m^TN_0V)^T\wedge(\om\times N_0)+dm^T\wedge
d(V+x_0)]+[\pa_w(m^TN_0V)^T(\om\times
N_0)-\pa_wm^Td(V+x_0)+dm^T\pa_wV]\wedge[m^TN_0V^T(\om\times
N_0)-m^Td(V+x_0)]=^{(\ref{eq:omom})}(N_0^T[\pa_w(m^TN_0V)\times(m^TN_0V)]dN_0
+N_0\times[-m^T\pa_wVd(m^TN_0V)+m^Td(V+x_0)\pa_w(m^TN_0V)-[\pa_wm^Td(V+x_0)
-dm^T\pa_wV]m^TN_0V])^T\wedge\om+[(d(V+x_0)\times\pa_wV)^T\wedge(m\times
dm)+\frac{1}{2}(\pa_wm\times m)^T(d(V+x_0)\times\wedge
d(V+x_0))]$. The last part is just (\ref{eq:inteco}); if in the
first part we replace $\wedge\om$ with $\odot dx_0$, then in the
obtained quantity the coefficients of the (linearly appearing)
second fundamental form of $x_0$ cancel, so it depends only on the
linear element of $x_0$; one can replace the last equation of
(\ref{eq:dis}) with this first part in which we ignore the
(linearly appearing) second fundamental form:
\begin{eqnarray}\label{eq:disTC}
m^T\pa_wV\neq 0,\ 2[d(V+x_0)\times\pa_wV]^TN_0\wedge
N_0^T[m\times(dm-m^TN_0dN_0)]+\nonumber\\
(\pa_wm\times m)^TN_0N_0^T[d(V+x_0)\times\wedge d(V+x_0)]\nonumber\\
+[\pa_w(m\times V)\times(m\times V)+m^T\pa_wVm\times
V]^TN_0N_0^T(dN_0
\times\wedge dN_0)=0,\nonumber\\
(N_0\times[m^T\pa_wVd(m^TN_0V)-(d(V+x_0)+V^TdN_0N_0)^T[m\pa_w(m^TN_0V)-\nonumber\\
\pa_wmm^TN_0V]-(dm-m^TN_0dN_0)^T\pa_wVm^TN_0V])^T\odot dx_0=0,\
du\wedge dv\wedge dw\neq 0.
\end{eqnarray}
If further we have the symmetric TC $m^TV=0$, then we can take
$m=:V\times N_0+\mathbf{m}N_0$ and (\ref{eq:disTC}) becomes
\begin{eqnarray}\label{eq:dissymTC}
V\times\pa_wV\neq 0,\ \frac{2[\pa_wV\times
d(V+x_0)]^T\wedge(V\times dx_0)}{(\pa_wV\times
V)^T(dx_0\times\wedge
dx_0)}+(\mathbf{m}^2+|V|^2)K=0,\nonumber\\
d\mathbf{m}=-\pa_w\mathbf{m}\frac{N_0^T[V\times
d(V+x_0)]}{N_0^T(\pa_wV\times
V)}+\mathbf{m}\frac{N_0^T(\pa_wV\times dV)}{N_0^T(\pa_wV\times
V)},\ du\wedge dv\wedge dw\neq 0.
\end{eqnarray}
If $\mathbf{m}=0$, then $(m\times V)\times N_0=0$ and as we shall
see later the $3$-dimensional IRDF must be a $2$-dimensional IRDF
counted with the simple $\infty$ of $w$ such that $x_0,x_0+V$ are
the focal surfaces of a normal congruence, which contradicts
$V\times\pa_wV\neq 0$.

Thus $\mathbf{m}\neq 0$; excluding the case $x_0$ developable we
can take $\mathbf{m}^2$ from the first equation of
(\ref{eq:dissymTC}) and replace it into the second one; applying
the compatibility condition $d\wedge$ to this equation and using
the equation itself we get $0=d\wedge
d\mathbf{m}=d\wedge(-\pa_w\mathbf{m}\frac{N_0^T[V\times
d(V+x_0)]}{N_0^T(\pa_wV\times
V)}+\mathbf{m}\frac{N_0^T(\pa_wV\times dV)}{N_0^T(\pa_wV\times
V)})=^{(\ref{eq:che})}(-\mathbf{m}(\frac{N_0^T(\pa_wV\times\pa_w^2V)}{N_0^T(\pa_wV\times
V)}\frac{N_0^T(V\times dV)}{N_0^T(\pa_wV\times
V)}+\frac{N_0^T(\pa_wV\times\pa_wdV)}{N_0^T(\pa_wV\times V)})
+\pa_w\mathbf{m}\frac{N_0^T[\pa_wV\times
d(V+x_0)]}{N_0^T(\pa_wV\times V)})\wedge\frac{N_0^T(V\times
dx_0)}{N_0^T(\pa_wV\times
V)}+\\\frac{-\mathbf{m}\pa_wV^TV+\pa_w\mathbf{m}|V|^2}{N_0^T(\pa_wV\times
V)}\frac{1}{2}N_0^T(dN_0\times\wedge dN_0)$.

Applying $\pa_w$ to the first equation of (\ref{eq:dissymTC}) we
get $(\frac{N_0^T(\pa_wV\times\pa_w^2V)}{N_0^T(\pa_wV\times
V)}\frac{N_0^T(V\times dV)}{N_0^T(\pa_wV\times
V)}+\frac{N_0^T(\pa_wV\times\pa_wdV)}{N_0^T(\pa_wV\times V)})
\wedge\frac{N_0^T(V\times dx_0)}{N_0^T(\pa_wV\times
V)}+\frac{N_0^T(\pa_wV\times dV)}{N_0^T(\pa_wV\times
V)}\wedge\frac{N_0^T(\pa_wV\times dx_0)}{N_0^T(\pa_wV\times
V)}+\frac{\mathbf{m}\pa_w\mathbf{m}+\pa_wV^TV}{N_0^T(\pa_wV\times
V)}N_0^T(dN_0\times\wedge dN_0)=0$; thus the previous relation
becomes
\begin{eqnarray}\label{eq:Wcoco}
\frac{2(\pa_wV\times dV)^T\wedge(\pa_wV\times dx_0)}{(\pa_wV\times
V)^T(dx_0\times\wedge
dx_0)}+(\mathbf{m}\pa_w\mathbf{m}+V^T\pa_wV)K=0,\ du\wedge
dv\wedge dw\neq 0
\end{eqnarray}
and (\ref{eq:dissymTC}) becomes
\begin{eqnarray}\label{eq:dissymTC1}
\pa_w\frac{N_0^T[\pa_wV\times d(V+x_0)]}{N_0^T(\pa_wV\times
V)}\wedge N_0^T(V\times dx_0)-\frac{N_0^T(\pa_wV\times
dV)}{N_0^T(\pa_wV\times V)}\wedge N_0^T(\pa_wV\times
dx_0)=0,\nonumber\\
\frac{1}{2}d(\frac{[2\pa_wV\times d(V+x_0)]^T\wedge(V\times
dx_0)}{K(\pa_wV\times V)^T(dx_0\times\wedge dx_0)}+|V|^2)=
-(\frac{2(\pa_wV\times dV)^T\wedge(\pa_wV\times
dx_0)}{K(\pa_wV\times V)^T(dx_0\times\wedge
dx_0)}+V^T\pa_wV)\nonumber\\\frac{N_0^T[V\times
d(V+x_0)]}{N_0^T(\pa_wV\times V)}+(\frac{2[\pa_wV\times
d(V+x_0)]^T\wedge(V\times dx_0)}{K(\pa_wV\times
V)^T(dx_0\times\wedge dx_0)}+|V|^2)\frac{N_0^T(\pa_wV\times
dV)}{N_0^T(\pa_wV\times V)},\nonumber\\ du\wedge dv\wedge dw\neq
0;
\end{eqnarray}
applying $\pa_w$ to the second equation of (\ref{eq:dissymTC1})
and using the first one we get another second order equation in
$V$.

\subsection{The Weingarten congruence property}

We consider a $3$-dimensional tangential IRDF with the symmetry of
the TC (thus $m=V\times N_0+\mathbf{m}N_0$ and $V,\mathbf{m}$
satisfy (\ref{eq:dissymTC})) and inquire in what case any
deformation $(x,dx)=(Rx_0+t,Rdx_0)$ of $x_0$ and any leaf are the
focal surfaces of a W congruence.

By the Darboux-Guichard's this is equivalent to the requirement
that for any leaf there is an infinitesimal deformation of $x$ in
the direction $Rm$ normal to the leaf, that is
$0=\frac{1}{2\rho}\frac{d}{d\ep}|_{\ep=0}|R^{-1}\ti d(x+\ep\rho
Rm)|^2=dx_0^T\odot(d\log\rho m+\om\times m+dm-\frac{m^T[\om\times
V+d(V+x_0)]}{m^T\pa_wV}\pa_wm)=dx_0^T(V\times N_0)\odot
(d\log\rho+\frac{N_0^T[\pa_wV\times d(V+x_0)]}{N_0^T(\pa_wV\times
V)}+\mathbf{m}\frac{\pa_wV^T(\om\times
N_0+dN_0)}{N_0^T(\pa_wV\times V)})$, or
\begin{eqnarray}\label{eq:Wco}
d\log\rho=-\frac{N_0^T[\pa_wV\times d(V+x_0)]}{N_0^T(\pa_wV\times
V)}-\mathbf{m}\frac{\pa_wV^T(\om\times
N_0+dN_0)}{N_0^T(\pa_wV\times V)}.
\end{eqnarray}
Imposing the compatibility condition $\ti d\wedge$ on
(\ref{eq:Wco}), using $dw=\frac{N_0^T[V\times
d(V+x_0)]}{N_0^T(\pa_wV\times V)}+\mathbf{m}\frac{V^T(\om\times
N_0+dN_0)}{N_0^T(\pa_wV\times V)}$ (note that for $w=w(u,v,c),\
c=$ct we have $d\log\pa_cw=-d\log\rho+...$, so up to a certain
scaling $\rho^{-1}$ can be interpreted as an infinitesimal
deformation of $w$), $(\om\times N_0+dN_0)\times\wedge(\om\times
N_0+dN_0)=dN_0\times\wedge dN_0$ and replacing $d\mathbf{m}$ from
the second equation of (\ref{eq:dissymTC}) we get
$0=-d\wedge\frac{N_0^T[\pa_wV\times d(V+x_0)]}{N_0^T(\pa_wV\times
V)}+\pa_w\frac{N_0^T[\pa_wV\times d(V+x_0)]}{N_0^T(\pa_wV\times
V)}\wedge\frac{N_0^T[V\times d(V+x_0)]}{N_0^T(\pa_wV\times
V)}+\mathbf{m}[\pa_w\frac{\mathbf{m}\pa_wV}{N_0^T(\pa_wV\times
V)}\times\frac{V}{N_0^T(\pa_wV\times
V)}]^T\frac{1}{2}(dN_0\times\wedge
dN_0)\\-\mathbf{m}([-\pa_w\log\mathbf{m}\frac{N_0^T[V\times
d(V+x_0)]}{N_0^T(\pa_wV\times V)}+\frac{N_0^T(\pa_wV\times
dV)}{N_0^T(\pa_wV\times V)}]\frac{\pa_wV}{N_0^T(\pa_wV\times
V)}+d\frac{\pa_wV}{N_0^T(\pa_wV\times
V)}-\pa_w\frac{N_0^T[\pa_wV\times d(V+x_0)]}{N_0^T(\pa_wV\times
V)}\\\frac{V}{N_0^T(\pa_wV\times
V)}+\frac{1}{\mathbf{m}}\frac{N_0^T[V\times
d(V+x_0)]}{N_0^T(\pa_wV\times
V)}\pa_w\frac{\mathbf{m}\pa_wV}{N_0^T(\pa_wV\times
V)})^T\wedge(\om\times N_0+dN_0)$. Using
$\pa_w\frac{N_0^T[\pa_wV\times d(V+x_0)]}{N_0^T(\pa_wV\times
V)}=\frac{N_0^T(\pa_wV\times\pa_wdV)}{N_0^T(\pa_wV\times
V)}+\frac{N_0^T[V\times d(V+x_0)]}{N_0^T(\pa_wV\times
V)}\frac{N_0^T(\pa_wV\times\pa_w^2V)}{N_0^T(\pa_wV\times V)}$ the
last term cancels and the remaining first three terms boil down to
(\ref{eq:Wcoco}).

Note also that the IC (\ref{eq:dissymTC}) can be easier obtained
by applying the compatibility condition $\ti d\wedge$ to
$dw=\frac{N_0^T[V\times d(V+x_0)]}{N_0^T(\pa_wV\times
V)}+\mathbf{m}\frac{V^T(\om\times N_0+dN_0)}{N_0^T(\pa_wV\times
V)}$ and using the equation itself: $0=\ti
d\wedge(\frac{N_0^T[V\times d(V+x_0)]}{N_0^T(\pa_wV\times
V)}+\mathbf{m}\frac{V^T(\om\times N_0+dN_0)}{N_0^T(\pa_wV\times
V)})=-\frac{\pa_w[N_0^T[V\times d(V+x_0)]+\mathbf{m}V^T(\om\times
N_0+dN_0)]}{N_0^T(\pa_wV\times V)}\wedge(\frac{N_0^T[V\times
d(V+x_0)]}{N_0^T(\pa_wV\times V)}+\mathbf{m}\frac{V^T(\om\times
N_0+dN_0)}{N_0^T(\pa_wV\times V)})+d\wedge\frac{N_0^T[V\times
d(V+x_0)]}{N_0^T(\pa_wV\times V)}+d\frac{mV^T}{N_0^T(\pa_wV\times
V)}\wedge(\om\times
N_0+dN_0)=^{(\ref{eq:che})}-[\frac{N_0^T[\pa_wV\times
d(V+x_0)]\wedge N_0^T(V\times dx_0)}{N_0^T(\pa_wV\times
V)}+(\mathbf{m}^2+|V|^2)\frac{1}{2}N_0^T(dN_0\times\wedge
dN_0)]\frac{1}{N_0^T(\pa_wV\times
V)}+(d\mathbf{m}+\pa_w\mathbf{m}\frac{N_0^T[V\times
d(V+x_0)]}{N_0^T(\pa_wV\times
V)}-\mathbf{m}\frac{N_0^T(\pa_wV\times dV)}{N_0^T(\pa_wV\times
V)})\wedge\frac{V^T(\om\times N_0+dN_0)}{N_0^T(\pa_wV\times V)}$.

\subsection{Applicability correspondence of leaves of a general nature}

We consider the question wether the leaves are deformations of
surfaces; we exclude the case $m$ isotropic because all leaves are
isotropic developables with degenerate $2$-dimensional linear
element.

The leaves for a particular deformation $x_0$ are applicable to
the same surface or to an $1$-dimensional family of surfaces.

Since we already have a correspondence between the facets of
leaves given by rolling, it is natural to require that the
applicability correspondence is independent of the shape of $x$.

However, for $(m\times V)\times N_0\neq 0$ the distribution of
facets into leaves changes with the shape of $x$; thus all leaves
of all rolled distributions of facets must be applicable to the
same surface (and we have a submersion from the $3$-dimensional
IRDF to the distribution of tangent planes of the fixed surface),
or $(m\times V)\times N_0=0$ and we have to discuss only a
$2$-dimensional IRDF with the leaf having linear element
independent of the shape of $x$.

Consider first the general case $m^T\pa_wV(m\times V)\times
N_0\neq 0$ and all leaves of all rolled distributions of facets
are applicable to the same surface
$y=y(u_1,v_1)=y(u_1(u,v,w),v_1(u,v,w))$; if for a particular
deformation of $x_0$ one knows the applicability correspondence of
all leaves to the surface $y$, then one finds the applicability
correspondence to $y$ of all rolled leaves (including the case of
leaves with degenerate linear element) by composing the rolling of
the particular leaves on $y$ with the inverse of the rolling of
$x_0$ on $x$. For $u_1,v_1=$ct from the three independent
variables $u,v,w$ only one remains independent, thus giving the
submersion from the IRDF to the distribution of tangent planes of
$y$ (equivalently we count each facet $(y,dy)$ with simple
$\infty$ multiplicity).

The function $w=w(u,v,c)$ is given by the integration of the
rolled (\ref{eq:inte}) $m^T(\om\times V+d(V+x_0)+\pa_wVdw)=0$;
thus a-priori $w$ depends also on $\om$ and we have
\begin{eqnarray}\label{eq:aplco}
|(I_3-\frac{\pa_wVm^T}{m^T\pa_wV})[\om\times V+d(V+x_0)]|^2=
|dy-\frac{m^T[\om\times
V+d(V+x_0)]}{m^T\pa_wV}\pa_wy|^2,\nonumber\\ \forall\ \om\
\mathrm{satisfying}\ (\ref{eq:om}).
\end{eqnarray}
The leaves are applicable to different regions of $y$ because the
constant $c$ in $w=w(u,v,c)$ changes for $\om$ fixed and for
different $\om$ $w$ changes; however in order to be determined by
the rolled (\ref{eq:inte}), $w$ is not allowed to be linked to
$\om$ by any other relation, either functional (as a-priori
(\ref{eq:aplco}) is) or differential; thus in (\ref{eq:aplco})
$\om$ cancels independently of $w$ and outside $w$ we can replace
$\om$ with any other $\hat\om$ satisfying (\ref{eq:om}). With
$\mathcal{M}:=(I_3-\frac{\pa_wVm^T}{m^T\pa_wV})^T(I_3-\frac{\pa_wVm^T}{m^T\pa_wV})
-\frac{|\pa_wy|^2}{(m^T\pa_wV)^2}mm^T,\
\mathcal{N}:=\frac{\pa_wy^Tdy}{m^T\pa_wV}m$ we have
$d(V+x_0)^T\odot[\mathcal{M}d(V+x_0)+2\mathcal{N}]
-|dy|^2=(\om\times V)^T\odot[\mathcal{M}(\om\times
V+2d(V+x_0))+2\mathcal{N}]=0,\ \forall\om$ satisfying
(\ref{eq:om}); in particular it is true if we replace $\om$ with
$\hat\om:=-2N_0\times dN_0,-\om-2N_0\times dN_0$; with
$\Om_0:=(N_0\times dN_0)\times V$ we have
$\Om_0^T\odot[\mathcal{M}(\Om_0-d(V+x_0))-\mathcal{N}]=(\om\times
V)^T\odot[\mathcal{M}(\Om_0-d(V+x_0))-\mathcal{N}]=(\om\times
V)^T\odot\mathcal{M}(\om\times V+2\Om_0)=0,\ \forall\om$
satisfying (\ref{eq:om}). For $\om+\ep\del\om$ infinitesimal
deformation of $\om$ (that is $(\del\om)^\bot=0,\
d\wedge\del\om+\del\om\times\wedge\om=0,\ \del\om\times\wedge
dx_0=0$) from the last relation we get $(\del\om\times
V)^T\odot\mathcal{M}(\om\times V+\Om_0)=0,\ \forall\om$ satisfying
(\ref{eq:om}) and $\del\om$ infinitesimal deformation of $\om$.

Thus excluding $x_0$ developable we need
$d(V+x_0)^T\odot[\mathcal{M}d(V+x_0)+2\mathcal{N}]-|dy|^2=0,\
[V\times \mathcal{M}(dx_0\times V)]\times N_0=[V\times
(\mathcal{M}d(V+x_0)+\mathcal{N})]\times N_0=0$.

If $V^TN_0\neq 0$, then we need $V\times\mathcal{M}(dx_0\times
V)=V\times[\mathcal{M}d(V+x_0)+\mathcal{N}]=0$; in particular we
get $0=(m\times V)^T\mathcal{M}(dx_0\times V)=(m\times
V)^T(I_3-\frac{\pa_wVm^T}{m^T\pa_wV})(dx_0\times V)$, so
$(I_3-\frac{\pa_wVm^T}{m^T\pa_wV})(dx_0\times V)\times\wedge
(I_3-\frac{\pa_wVm^T}{m^T\pa_wV})(dx_0\times V)=(dx_0\times\wedge
dx_0)^TN_0N_0^TV[(I_3-\frac{\pa_wVm^T}{m^T\pa_wV})^*]^TV$ is a
multiple of $m\times V$; since
$\ker[(I_3-\frac{\pa_wVm^T}{m^T\pa_wV})^*]^T=
\mathrm{Im}(I_3-\frac{\pa_wVm^T}{m^T\pa_wV})^T=(\pa_wV)^\bot,\
\mathrm{Im}[(I_3-\frac{\pa_wVm^T}{m^T\pa_wV})^*]^T=\mathbb{C}m$
and excluding $m$ isotropic we get $V^T\pa_wV=0$ and further
$\pa_wV=\mathbf{v}\times V,\ \mathbf{v}^TN_0=0$; similarly we get
$\pa_wV^T[d(V+x_0)\times\wedge d(V+x_0)]=0$, that is the leaves
are (isotropic) curves or points. But the choice of surface $x_0$
was irrelevant, so this property must be true if we replace $x_0$
with $x$ and $V$ with $RV$, that is $0=\pa_wV^T[(\om\times
V+d(V+x_0))\times\wedge (\om\times
V+d(V+x_0))]=-2N_0^TVN_0^T(\om\times\mathbf{v})\wedge
V^Td(V+x_0),\ \forall\ \om$ satisfying (\ref{eq:om}), or
$N_0^T(dx_0\times\mathbf{v})\wedge V^Td(V+x_0)=0$; such a
configuration, even if it existed, does not apply to our
considerations since we are looking for leaves which are surfaces.

Thus we have the TC $V^TN_0=0$ and we need
$d(V+x_0)^T\odot[\mathcal{M}d(V+x_0)+2\mathcal{N}] -|dy|^2=0,\
N_0^T\mathcal{M}N_0=0,\ N_0^T[\mathcal{M}d(V+x_0)+\mathcal{N}]=0$,
so $|\pa_wy|^2=\pa_wV^T(I_3+\frac{mm^T}{(m^TN_0)^2})\pa_wV,\
\pa_wy^Tdy=d(V+x_0)^T(I_3-N_0N_0^T)(I_3+\frac{mm^T}{(m^TN_0)^2})\pa_wV,\
|dy|^2=d(V+x_0)^T(I_3-N_0N_0^T)(I_3+\frac{mm^T}{(m^TN_0)^2})
(I_3-N_0N_0^T)d(V+x_0)$ (note $m^TN_0=0\Rightarrow (m\times
V)\times N_0=0$); thus $|\ti dy|^2=\ti
d(V+x_0)^T(I_3-N_0N_0^T)(I_3+\frac{mm^T}{(m^TN_0)^2})
(I_3-N_0N_0^T)\ti d(V+x_0)=|(I_3-\frac{N_0m^T}{m^TN_0})\ti
d(V+x_0)|^2$ and
\begin{eqnarray}\label{eq:tidy}
\ti dy=R_1(I_3-\frac{N_0m^T}{m^TN_0})\ti d(V+x_0),\
R_1\subset\mathbf{O}_3(\mathbb{C}),\ du\wedge dv\wedge dw\neq 0.
\end{eqnarray}
Note that since along the leaves we have $m^T\ti d(V+x_0)=0$,
$R_1$ is the rotation of the rolling of the leaves on $y$; if we
replace $x_0$ with an applicable surface $x=Rx_0+t$, then $R_1$ is
replaced with $R_1R^{-1}$.

Imposing the compatibility condition $R_1^{-1}\ti d\wedge$ on
(\ref{eq:tidy}) we get $0=[R_1^{-1}\ti
dR_1(I_3-\frac{N_0m^T}{m^TN_0})-\ti
d(\frac{N_0m^T}{m^TN_0})]\wedge\ti d(V+x_0)$; with $R_1^{-1}\ti
dR_1=:\Om_1du+\Om_2dv+\Om_3dw$ this constitutes a linear system of
$9$ equations on the $9$ entries of $\Om_j,\ j=1,2,3$ with the
rank of the matrix of the system being $6$, so the rank of the
augmented matrix of the system must be also $6$ and the solution
$R_1^{-1}\ti dR_1$ must also satisfy the compatibility condition
$\ti d\wedge(R_1^{-1}\ti dR_1)+\frac{1}{2}[R_1^{-1}\ti dR_1,\wedge
R_1^{-1}\ti dR_1]=0$; these are the necessary and sufficient
conditions on the $3$-dimensional tangential IRDF in order to
obtain applicability correspondence of leaves of a general nature.

By applying, if necessary, a change of variable $w=w(\ti w,u,v)$,
we have $N_0^T[d(V+x_0)\times\wedge d(V+x_0)]\neq 0$ and the above
considered linear system is consistent for
\begin{eqnarray}\label{eq:cons}
d(V+x_0)^T(I_3-\frac{mN_0^T}{m^TN_0})\odot
[d(\frac{N_0m^T}{m^TN_0})\pa_wV-\pa_w(\frac{N_0m^T}{m^TN_0})d(V+x_0)+\nonumber\\
2N_0^T[\pa_wV\times d(V+x_0)]\frac{d(\frac{N_0m^T}{m^TN_0})\wedge
d(V+x_0)}{N_0^T[d(V+x_0)\times\wedge d(V+x_0)]}]=0,\ du\wedge
dv\wedge dw\neq 0;
\end{eqnarray}
if we further assume the symmetric TC $m^TV=0$, then by using
(\ref{eq:dissymTC}) (\ref{eq:cons}) is satisfied.

\subsection{The singular cases}

The law (\ref{eq:inte}) of distribution of facets into leaves is
independent of rolling if $0=m^T(\om\times V+d(V+x_0)+\pa_wVdw)$
is independent of $\om,\ \forall\om$ satisfying (\ref{eq:om}),
that is for $(m\times V)\times N_0=0\Leftrightarrow
V^TN_0=m^TN_0=0\vee m=V$, in which case we have an arbitrary
$1$-dimensional family of $2$-dimensional IRDF ($w$ can be
prescribed in any continuous manner).

In the case of $2$-dimensional IRDF we have $0=m^T(\om\times
V+d(V+x_0)),\ \forall\om$ satisfying (\ref{eq:om}), so $(m\times
V)\times N_0=0$ and for the cancelling of the coefficients of the
linearly appearing second fundamental form we get the vacuous
$0=[d(m\times V)-(N_0\times dN_0)\times(m\times
V)]^T\wedge(N_0\times dN_0)$; thus we need only
\begin{eqnarray}\label{eq:dis1}
(m\times V)\times N_0=0,\ 2(dm-(N_0\times dN_0)\times
m)^T\wedge(d(V+x_0)-(N_0\times dN_0)\times V)+\nonumber\\
(m\times V)^T(dN_0\times\wedge dN_0)=0.
\end{eqnarray}
In the case $m=V$ we have $V^Td(x_0+V)=0$, so $|V|^2\neq 0$ and
$x_0+sV$ forms a normal congruence; if $V^TN_0\neq 0$, then we
have envelopes of sphere congruences (facets are tangent to
spheres centered on $x_0$).

Dupin and Malus studied normal congruences; they remain normal
after reflections and refractions in surfaces and this property is
independent of the shape of the surface (that is if we
transversely capture a normal congruence in $x_0$, deform $x_0$
and release the congruence after a constant angle refraction law,
it remains normal), which explains the fact that envelopes of
sphere congruences centered on a surface are independent of the
shape of the surface (Beltrami); conversely Levi-Civita proved
that any two normal congruences can be transformed one into the
other by two reflections or refractions in surfaces.

If we tangentially capture a normal congruence in $x_0$ ($m=V,\
V^TN_0=0$), then $x_0$ is a focal surface of the normal
congruence. The developables of the normal congruence are
generated by varying the normals on a normal surface along the
lines of curvature; they envelope an $1$-dimensional family of
curves on each of the focal surfaces; thus this system of curves
give a conjugate system on both focal surfaces. Since the tangent
planes of the two focal surfaces are generated by the normals of a
normal surface and the tangents of one of its lines of curvature,
the osculating plane of a curve on a focal surface (whose tangent
surface is one of the developables of the normal congruence) is
normal to the tangent plane of the respective focal surface and
thus the curve is a geodesic (see Eisenhart (\cite{E1},\S 74)).

If we deform $x_0$ and release the congruence, then it remains
normal with the same focal surfaces (but the linear element of the
other focal surface and the other curves of the conjugate system
change except for normal W congruences (see Bianchi (\cite{B5},Vol
{\bf 4},(202)))): take an $1$-dimensional family of geodesics
$v=$ct on $x_0$ and their orthogonal trajectories $u=$ct; the
normal congruence is formed by the tangents to geodesics. We have
the linear element $|dx_0|^2=du^2+G(u,v)dv^2$ such that the curves
$v=$ct are geodesics on $x_0$. From $p=x_0+s(u,v)\pa_ux_0,\
0=\pa_ux_0^T(dp\wedge\times
dp)=s(\pa_ux_0\times\pa_u^2x_0)^T(\pa_vx_0+s\pa_{uv}^2x_0)du\wedge
dv$ we get $p=x_0-\frac{(\pa_ux_0\times\pa_u^2x_0)^T\pa_vx_0}
{(\pa_ux_0\times\pa_u^2x_0)^T\pa_{uv}^2x_0}\pa_ux_0=x_0-2\frac{G}{\pa_uG}\pa_ux_0$
(we exclude the case of developables $\pa_uG=0$, when such a
surface $p$ does not exist); the normal surfaces are given by
$x_0+(c-u)\pa_ux_0,\ c=$ct, so $u$ is, up to addition with a
constant, a principal curvature (see Eisenhart (\cite{E1},\S 76)).
In particular since by Chasles's theorem the common tangents to
$2$ confocal quadrics form a normal congruence that envelopes
geodesics on the $2$ confocal quadrics, if we capture this normal
congruence in a quadric, deform the quadric and release the
congruence, then it remains normal with the same focal surfaces
(but the linear element of the other focal surface changes). The
linear element of the other focal surface is independent of the
shape of $x$ only for normal W congruences, when both focal
surfaces are applicable to surfaces of revolution, with the
geodesics on $x_0$ corresponding to meridians.

In the case $V^TN_0=m^TN_0=0,\ |m|^2\neq 0$ the configuration is
obtained by taking $x_0$ the envelope of a family of normal planes
to the leaf $p$ (planes containing the direction $m$; to each
point of $p$ corresponds a plane); the facets and the surface
$x_0$ are thus individuated and if we deform $x_0$, then the
rolled distribution of facets is still integrable. If the
particular leaf $p$ is a curve, then we take the intersection of
tangent planes of $x_0$ with $p$ and the facet is tangent to $p$
and normal to the tangent plane of $x_0$; if the particular leaf
$p$ is a point, then $x_0$ must be a cone passing through $p$.

If one of the isotropic directions in each facet of the
distribution can be brought for a particular deformation $x_0$ to
coincide for all points of $x_0$, then the leaf in this particular
position is an isotropic line (since it cannot be isotropic
developable). In this case we can take any isotropic developable
containing this isotropic line, intersect the tangent planes of
$x_0$ with it to obtain curves $c=c(u,v,w)$ in $Tx_0$ such that
the $3$-dimensional distribution formed by the facets normal
planes of $c$ is integrable for all deformations $x$ of $x_0$
(Ribaucour and Darboux; see Darboux (\cite{D1},\S\ 762)). Such is
the case for {\it cyclic systems} (the generating isotropic
developable is a null cone and the curves are circles in tangent
planes of $x$). If $m^TV=0$, then $V$ is a normal congruence
(captured tangentially in $x_0$) with $V+x_0$ being the other
focal surface; if further this is a normal W congruence, then the
construction above is possible and we get the deformation of a
surface of revolution to an isotropic line (see Darboux
(\cite{D2},\S 169)).

In what concerns ACLGN assume that we have all rolled leaves
applicable without collapsing ansatz:
\begin{eqnarray}\label{eq:aplco1}
(\om\times V+d(V+x_0))^T(\om\times
V+d(V+x_0))=d(V+x_0)^Td(V+x_0)\ \mathrm{non-degenerate},\nonumber\\
\forall\om\ \mathrm{satisfying\ (\ref{eq:om})}.
\end{eqnarray}
This becomes $(\om\times V)^T\odot(\om\times V+2d(V+x_0))=0,\
\forall\om$ satisfying (\ref{eq:om}); in particular it is true for
$\om:=-2N_0\times dN_0,-\om-2N_0\times dN_0$; with
$\Om_0:=(N_0\times dN_0)\times V$ we have
$\Om_0^T\odot(\Om_0-d(V+x_0))=0,\ (\om\times V)^T\odot(\om\times
V+2\Om_0)=0,\ \forall\om$ satisfying (\ref{eq:om}). For
$\om+\ep\del\om$ infinitesimal deformation of $\om$ (that is
$(\del\om)^\bot=0,\ d\wedge\del\om+\del\om\times\wedge\om=0,\
\del\om\times\wedge dx_0=0$) we get $(\del\om\times
V)^T\odot(\om\times V+\Om_0)=0,\ \forall\om$ satisfying
(\ref{eq:om}) and $\del\om$ infinitesimal deformation of $\om$, a
contradiction.

Thus (\ref{eq:aplco1}) must be refined to assume collapsing ansatz
of leaves: for a particular deformation of $x_0$ (which can be
taken to be $x_0$) the leaf has degenerate $2$-dimensional element
and must be a(n isotropic) curve or a point.

For $m=V$ the sphere congruence is defined so that the spheres
centered on $x_0$ are tangent to the given (isotropic) curve or
point; the linear element of the leaves depends however on the
shape of $x$.

Note however that envelopes of sphere congruences enjoy other
properties; it is worth mentioning here the sequence of events
that led Bianchi to his discovery of the B transformation of
quadrics: in 1899 (and based on an earlier result of 1897)
Guichard discovered that when a quadric with(out) center and of
revolution around the focal axis rolls on one of its deformations,
its foci (focus) describe CMC (minimal) surfaces (thus
generalizing an earlier result of Bonnet on the rolling of the
unit sphere on a CGC $1$ surface; since a surface of revolution
can be rolled on the axis of revolution with the arc-length of the
meridian corresponding to the arc-length of the axis, this is also
a generalization of Delaunay's generation of CMC (minimal)
surfaces of revolution with the meridian being described by a
focus of a conic as the conic rolls on the axis in a meridian
plane) and the same result for the intersections of the isotropic
rulings of the quadric with the tangent planes of the quadric (in
this case according to Darboux the CMC (minimal) surfaces are the
rolled foci (focus) as the quadric rolls on the complementary
transform of the considered deformation).

Building on Guichard's result Darboux reduced the deformations of
the Darboux quadrics with center to the deformations of the
(pseudo-)sphere (Goursat had integrated earlier the equations for
the deformations of certain Darboux paraboloids): if we intersect
the tangent planes of a surface with a(n isotropic) plane and
consider the resulting congruence of lines as the surface rolls on
one of its deformations, then the focal surfaces of this
congruence are obtained from the intersection of the common
conjugate directions with the lines and the conjugate system
induced by the developables of the congruence on the focal
surfaces corresponds to the conjugate system common to the surface
and its deformation; if the plane is isotropic, then the
congruence is normal and the parallel surfaces are given by the
intersection of the isotropic lines of the isotropic plane with
the line of the congruence (thus they are envelopes of rolling
congruences of the isotropic lines of the isotropic plane).

In Guichard's result with $4$ isotropic rulings (situated in two
isotropic planes) the $4$ points in the tangent plane are situated
on the rulings at the tangency point; since the asymptotes and any
conjugate system are harmonically conjugate we get two normal
congruences each having two parallel surfaces in harmonic ratio
with the focal surfaces and this configuration is possible only if
the two parallel surfaces have CMC (this configuration is
preserved under conformal changes of the space that preserve
geodesics, so we get parallel CMC surfaces in space forms by
changing Cayley's absolute; see Bianchi (\cite{B5},Vol {\bf
4},(108))).

Thus Guichard's for isotropic rulings remains valid if we consider
only two parallel isotropic rulings in an isotropic plane, that is
a Darboux quadric with center.

Darboux inquired what becomes of Guichard's result if we consider
the intersection of the $8$ isotropic rulings on the general
quadric with tangent planes; the resulting surfaces (envelopes of
rolling congruences of isotropic rulings of the quadric) are {\it
isothermic} (surfaces with isothermal lines of curvature) with the
lines of curvature corresponding to the conjugate system common to
the quadric and its deformation (thus in conformal
correspondence); however, since the deformation problem depends on
two functions of a variable and the space of isothermic surfaces
in conformal correspondence is much larger, Darboux sought to find
the properties that individuate the special isothermic surfaces
obtained above.

The surfaces corresponding to parallel isotropic rulings have same
normal direction and form a harmonic ratio with the focal surfaces
of the congruence of their joins (thus giving the involutory
Christoffel transformation of isothermic surfaces) and the
surfaces corresponding to isotropic rulings that intersect (at
umbilics) are leaves of the cyclic system generated by the
isotropic cones at umbilics; thus they are envelopes of a sphere
congruence.

Conversely, Darboux proved that any isothermic surface appears as
envelope of a $4$-dimensional family of sphere congruences whose
other envelopes are isothermic surfaces in conformal
correspondence with the given isothermic surface; thus introducing
the {\it Darboux} (D) transformation of isothermic surfaces (a
particular case of {\it Ribacour transformation} (envelopes of
sphere congruences with correspondence of lines of curvature; note
however that this property is independent of the shape of the
surface of centers only if it is applicable to a quadric of
revolution), which in turn corresponds via Lie's contact
transformation (which exchanges the points and planes of facets)
to the B transformation of the focal surfaces of a W congruences);
later on he found the sought characterization of special
isothermic surfaces and the fact that in the $4$-dimensional space
of D transforms of a special isothermic surface a $3$-dimensional
sub-space are special isothermic associated to the same quadric;
thus the space of special isothermic surfaces associated to the
same quadric remains closed under the iteration of the D
transformation, which provided the first theory of deformations of
general quadrics depending on arbitrarily many constants.

Building on Guichard's and Darboux's results Bianchi proved in
1899 the inversion of Guichard's result (all CMC (minimal)
surfaces can be realized in a $2$-dimensional fashion as in
Guichard's result and if a point rigidly attached to a surface
generates CMC (minimal) surfaces when the surface rolls on its
deformations, then the surface and point must be as in Guichard's
result); in 1904 he proved the BPT for the D transformation of
isothermic surfaces (the iteration of the D transformation can be
realized using only algebraic and differential computations) and
provided a simple geometric interpretation of Darboux's
characterization of special isothermic surfaces: they are the
closure of CMC surfaces under conformal transformations of the
space (see Bianchi (\cite{B5},Vol {\bf 4},(108))).

In trying to see what becomes of Darboux's result for CMC
$\frac{1}{2}$ surfaces he realized that the transformation induced
for CGC $1$ surfaces is not a fundamental one, but rather the
composition of two such ones; thus he found in 1899 the B
transformation of the sphere (the reason this B transformation was
discovered much later than that of the pseudo-sphere is that at
the first iteration it gives complex leaves; one needs the
iteration of two complex conjugate B transformations (via the BPT)
to get back real surfaces), which he later generalized to quadrics
of revolution and paraboloids (1905); once he realized that the
applicability law of the leaves is given at the level of confocal
quadrics by the Ivory affinity he generalized the B transformation
to all quadrics \cite{B1}.

For $m^TN_0=V^TN_0=0$ consider the defining surface $x_0$ such
that the leaf is a curve $c=c(s)$; the tangent planes of $x_0$ cut
$c$ at points $p$ and at those points $m=N_0\times c'(s)$. We have
$0=N_0^T(c(s)-x_0)$, so $ds=-\frac{dN_0^TV}{N_0^Tc'(s)}$ and the
linear element $|\om\times
V-\frac{dN_0^TV}{N_0^Tc'(s)}c'(s)|^2=|V^T(\om\times N_0+dN_0)|^2
+\frac{(dN_0^TV)^2}{(N_0^Tc'(s))^2}|c'(s)\times N_0|^2$  of the
rolled leaf must be non-degenerate independent of $\om$ for most
$\om\neq 0,-2N_0\times dN_0$ satisfying (\ref{eq:om}). Using
(\ref{eq:secoi}) this becomes $((N_0\times V)^TR^{-1}dN)^2$ being
independent of the shape of $x$ for most deformations $x$ of
$x_0$; if we refer $x_0$ to a system of coordinates formed by the
curves $u=$ct envelopes of the tangent field $N_0\times V$ and
their orthogonal trajectories $v=$ct, then we need
$(N^Td\pa_vx)^2$ to be independent of shape of $x$ for most
deformations $x$ of $x_0$, so we get a normal W congruence.

If the original leaf is a point $p$, then $N_0^T(p-x_0)=0,\
dN_0^T(p-x_0)=0$ and $x_0$ is developable with rulings passing
through $p$.

For $(m\times V)\times N_0\neq 0,\ m^T\pa_wV=0$ to each point of
$x_0$ we associate an $1$-dimensional family of facets centered on
a curve (or at a point) in a leaf. If the leaves of the original
distribution are an $1$-dimensional family of points, curves or
surfaces, then we can take the parametrization $(u,v)$ on $x_0$
such that as $u$ varies the curve of facets varies in its leaf and
as $v$ varies the curve of facets varies in different leaves; if
the original distribution is just the facets tangent to a curve,
each facet being counted with simple $\infty$ and such that the
facets associated to a point of $x_0$ are the simple $\infty$ of
facets centered at a point (counted with the simple $\infty$ of
$w$), then again we can take the parametrization $(u,v)$ on $x_0$
such that when $u$ varies we get the simple $\infty$ multiplicity
of the facet and when $v$ varies we get facets centered at
different points; the remaining case is when the original
distribution is just the facets tangent to a curve or surface,
each facet being counted with simple $\infty$ multiplicity, but
such that the facets associated to a point of $x_0$ are centered
on a curve.

However, only the case of the symmetric TC is pertinent to our
considerations and we shall not discuss it here.

\section{$3$-dimensional integrable tangential rolling
distributions of facets with collapsing ansatz}

In the case of collapsing ansatz of leaves (for a particular
deformation of $x_0$ (which can be taken to be $x_0$) the leaves
of the $3$-dimensional IRDF collapse to an $1$-dimensional family
of points or curves) without the TC one can derive the IC by
requiring that the reflected original distribution is also
integrable with the cancelling of the (linearly appearing) terms
containing the second fundamental form of $x_0$.

However, since we are interested in ACLGN this IC does not apply
to our considerations, so we have to further assume the TC.

Assume that we have a $3$-dimensional tangential distribution of
facets
$$(p,P)=(p(u_0,v_0,w),P(u_0,v_0,w)),\ p\in P$$ distributed in
the tangent planes of the surface $x_0=x_0(u_0,v_0):\ V:=p-x_0, \
V^TN_0=0$.

Further assume that the IC of this distribution depends only on
the linear element of $x_0$, that is if we roll $x_0^0\subseteq
x_0$ on one if its deformations
$(x^0,dx^0)=(R_0x_0^0+t_0,dx_0^0)$, then the rolled distribution
of facets $(R_0p+t_0,RP)=(R_0V+x^0,R_0P)$ remains integrable with
leaves $x^1$ (and normal fields $R_0m_0^1$).

Further assume that for one of the rolled distributions of facets
$(R_0p+t_0,R_0P)$ the leaves collapse from surfaces to curves (if
they collapse to points, then the rolled leaves are always
curves); we can take that defining distribution to be the initial
one and thus $x_0$ is the defining surface of the IRDF.

The leaves of the defining distribution of facets are an
$1$-dimensional family of curves that generate an auxiliary
surface $x_z$ (assume that we don't have a curve counted with
simple $\infty$ multiplicity, in which case the rolled leaves are
always curves).

If we reflect the distribution of facets in the tangent bundle of
$x_0$ (that is we roll $x_0$ on its other side), then we get a
$3$-dimensional integrable distribution of facets centered on
$x_z$, so the leaves are also curves that generate $x_z$.

Excluding the case $(N_0^0)^Tm_0^1=0$ (when the two families of
curves coincide and we get an $1$-dimensional family of
$2$-dimensional IRDF) and the case when $x_0=x_z$ is a ruled
surface with the leaves of the distribution being the rulings of
$x_0$ (note that we also have to exclude $x_0$ developable, in
which case by deforming it to a plane we get the rolled $x_z$
(captured in the $1$-dimensional family of tangent planes of
$x_0$) the same plane, so the plane itself is defining and
auxiliary surface, but the planes of the facets (which must be
tangent planes of curves in this plane) must leave the plane, so
we are in the previous situation), the two families of curves give
a parametrization $(u_1,v_1)$ of $x_z$ suited to our purposes;
$(V_0^1)^TN_0^0=0,\ V_0^1:=x_z^1-x_0^0$ imposes a functional
relationship between the four independent variables
$u_0,v_0,u_1,v_1$, leaving only three of them independent and any
other needed functional relationship between the four independent
variables $u_0,v_0,u_1,v_1$ must be a consequence of
$(V_0^1)^TN_0^0=0$

We have the reflection property
\begin{eqnarray}\label{eq:mr}
(m_0^1)^T\pa_{u^1}x_z^1=(m_0^1)^T(I_3-2N_0^0(N_0^0)^T)\pa_{v^1}x_z^1=0
\end{eqnarray}
and $R_0^{-1}\ti dx^1=\ti d(V_0^1+x_0^0)+R_0^{-1}dR_0V_0^1=\ti
dx_z^1+\om_0\times V_0^1,\ \om_0:=N_0^0\times R_0^{-1}dR_0N_0^0$.
But $(\om_0)^{\bot}=0$ and $\ti
dx_z^1=\pa_{u_1}x_z^1du_1+\pa_{v_1}x_{z}^1dv_1$, so
$0=(R_0m_0^1)^Tdx^1$ (and similarly $0=(R_0'{m_0^1}')^Tdx^1$ for
$R_0'=R_0(I_3-2N_0^0(N_0^0)^T),\ \om_0'=-\om_0-2N_0^0\times
dN_0^0$) becomes:
\begin{eqnarray}\label{eq:dify1}
-(V_0^1)^T(\om_0\times
N_0^0)(m_0^1)^TN_0^0+(m_0^1)^T\pa_{v_1}x_z^1dv_1=0\nonumber\\
(\Leftrightarrow -(V_0^1)^T(\om_0'\times
N_0^0)({m_0^1}')^TN_0^0+({m_0^1}')^T\pa_{u_1}x_z^1du_1=0)\Leftrightarrow\
B_z\ \mathrm{transformation}.
\end{eqnarray}
We get the transformation $B_z'$ for the reflected rolled
distribution of facets by the change
$(u_1,m_0^1)\leftrightarrow(v_1,{m_0^1}')$ or
$\om_0\leftrightarrow\om_0'$. Since the distributions of facets
$\mathcal{D},\mathcal{D'}$ with normal fields $m_0^1,{m_0^1}'$
reflect in $Tx_0^0$, the rolled distributions of facets
$(R_0,t_0)\mathcal{D},(R_0,t_0)\mathcal{D'}$ reflect in $Tx^0$, so
$B'_z(x^0)$ is just $B_z(x^0)$ when $x_0^0$ rolls on the other
face of $x^0$.

Using (\ref{eq:mr}) (\ref{eq:dify1}) becomes:
\begin{eqnarray}\label{eq:dify2}
\ \ \ \ \ \ -(V_0^1)^T(\om_0\times
N_0^0)+2(\pa_{v_1}x_z^1)^TN_0^0dv_1=0\ (\Leftrightarrow
-(V_0^1)^T(\om_0'\times N_0^0)+2(\pa_{u_1}x_z^1)^TN_0^0du_1=0).
\end{eqnarray}
(because $0=\ti
d((V_0^1)^TN_0^0)=(\pa_{v_1}x_z^1)^TN_0^0dv_1+(\pa_{u_1}x_z^1)^TN_0^0du_1
+(V_0^1)^TdN_0^0$ the equivalency of the equations of
(\ref{eq:dify2}) is clear).

Imposing the compatibility condition $\ti d\wedge$ on
(\ref{eq:dify2}) and using the equation itself we get the IC
$0=-(\pa_{u_1}x_z^1du_1+\pa_{v_1}x_z^1dv_1)^T\wedge(\om_0\times
N_0^0)+2(\pa_{u_1v_1}^2x_z^1)^TN_0^0du_1\wedge
dv_1+2(\pa_{v_1}x_z^1)^TdN_0^0\wedge
dv_1=^{(\ref{eq:omom})}\frac{1}{2}[V_0^1\times
(\frac{\pa_{u_1}x_z^1}{(\pa_{u_1}x_z^1)^TN_0^0}
+\frac{\pa_{v_1}x_z^1}{(\pa_{v_1}x_z^1)^TN_0^0})]^TN_0^0(dN_0^0)^T\wedge\om_0
+(dN_0^0)^TV_0^1[N_0^0\times(\frac{\pa_{u_1}x_z^1}{(\pa_{u_1}x_z^1)^TN_0^0}
+\frac{\pa_{v_1}x_z^1}{(\pa_{v_1}x_z^1)^TN_0^0})]^T\wedge\om_0-
(\pa_{u_1v_1}^2x_z^1)^TN_0^0\frac{(dN_0^0)^TV_0^1\wedge(N_0^0\times
V_0^1)^T\om_0}{(\pa_{u_1}x_z^1)^TN_0^0(\pa_{v_1}x_z^1)^TN_0^0}=
[\frac{1}{2}(V_0^1\times
N_0^0)^T(\frac{\pa_{u_1}x_z^1}{(\pa_{u_1}x_z^1)^TN_0^0}
+\frac{\pa_{v_1}x_z^1}{(\pa_{v_1}x_z^1)^TN_0^0})dN_0^0-\\
(dN_0^0)^T(\frac{\pa_{u_1}x_z^1}{(\pa_{u_1}x_z^1)^TN_0^0}
+\frac{\pa_{v_1}x_z^1}{(\pa_{v_1}x_z^1)^TN_0^0}-
\frac{(\pa_{u_1v_1}^2x_z^1)^TN_0^0V_0^1}
{(\pa_{u_1}x_z^1)^TN_0^0(\pa_{v_1}x_z^1)^TN_0^0})(V_0^1\times
N_0^0)]^T\wedge\om_0$, or $(V_0^1)^TN_0^0=0\Rightarrow
[(V_0^1\times N_0^0)^T(\pa_{u_1}x_z^1(\pa_{v_1}x_z^1)^T
+\pa_{v_1}x_z^1(\pa_{u_1}x_z^1)^T)N_0^0dN_0^0-
2(dN_0^0)^T(\pa_{u_1}x_z^1(\pa_{v_1}x_z^1)^T
+\pa_{v_1}x_z^1(\pa_{u_1}x_z^1)^T-V_0^1(\pa_{u_1v_1}^2x_z^1)^T)N_0^0(V_0^1\times
N_0^0)]^T\odot dx_0^0=0$. Referring $x_0^0$ to its lines of
curvature parametrization we get (excluding $x_0^0$ developable)
the IC
\begin{eqnarray}\label{eq:int}
(V_0^1)^TN_0^0=0\Rightarrow
(I_3-N_0^0(N_0^0)^T)[\pa_{u_1}x_z^1(\pa_{v_1}x_z^1)^T
+\pa_{v_1}x_z^1(\pa_{u_1}x_z^1)^T]N_0^0=(\pa_{u_1v_1}^2x_z^1)^TN_0^0V_0^1.\nonumber\\
\end{eqnarray}
If $\pa_{u_1v_1}^2x_z=0$, then $x_z(u_1,v_1)=f(u_1)+g(v_1),\
f'(u_1)\times g'(v_1)\neq 0$ and $(V_0^1)^TN_0^0=0\Rightarrow
f'(u_1)^TN_0^0=g'(v_1)^TN_0^0=0$ or
$f'(u_1)^TN_0^0g'(v_1)^TN_0^0\neq 0,\
\frac{f'(u_1)}{f'(u_1)^TN_0^0}+\frac{g'(v_1)}{g'(v_1)^TN_0^0}=2N_0^0,\
[f'(u_1)\times g'(v_1)]^TN_0^0=0$; in the first case we have
$N_0^0=\frac{f'(u_1)\times g'(v_1)}{|f'(u_1)\times g'(v_1)|}$ and
thus constant along the curve of tangency of the tangent cone of
$x_0$ from $x_z(u_1,v_1)$ and $x_0$ is developable and in the
second case $N_0^0$ along the curve of tangency of the tangent
cone of $x_0$ from $x_z(u_1,v_1)$ must span a line or a plane;
since the tangent cone of $x_0$ from $x_z(u_1,v_1)$ cannot be a
cylinder, $x_0$ is again developable.

Note also $(\pa_{u_1}x_z^1\times V_0^1)\times[\pa_{u_1}x_z^1\times
(I_3-2N_0^0(N_0^0)^T)\pa_{v_1}x_z^1]=-[(V_0^1\times
N_0^0)^T(\pa_{u_1}x_z^1(\pa_{v_1}x_z^1)^T
+\pa_{v_1}x_z^1(\pa_{u_1}x_z^1)^T)N_0^0
+(V_0^1)^TN_0^0(\pa_{u_1}x_z^1\times\pa_{v_1}x_z^1)^TN_0^0]\pa_{u_1}x_z^1$,
so we are actually in the symmetric TC case $(V_0^1)^Tm_0^1=0$; we
can take
\begin{eqnarray}\label{eq:mm'}
m_0^1:=\mathcal{B}_1\pa_{u_1}x_z^1\times V_0^1,\
{m_0^1}':=\mathcal{B}_1\pa_{v_1}x_z^1\times V_0^1,\
\mathcal{B}_1:=-z[(\pa_{u_1}x_z^1)^TN_0^0(\pa_{v_1}x_z^1)^TN_0^0]^{-1}.
\end{eqnarray}
Multiplying the first equation of (\ref{eq:dify2}) with
$\mathcal{B}_1(\pa_{u_1}x_z^1)^TN_0^0$ and using
$-\mathcal{B}_1(\pa_{u_1}x_z^1)^TN_0^0V_0^1=\\-\mathcal{B}_1(\pa_{u_1}x_z^1\times
V_0^1)\times N_0^0=-m_0^1\times N_0^0$ we get
\begin{eqnarray}\label{eq:BB'}
(m_0^1)^T\om_0+2zdv_1=0\ \Leftrightarrow
({m_0^1}')^T\om'_0+2zdu_1=0\ \Leftrightarrow B_z\
\mathrm{transformation},\nonumber\\
(m_0^1)^T\om_0'+2zdv_1=0\ \Leftrightarrow
({m_0^1}')^T\om_0+2zdu_1=0\ \Leftrightarrow B_z'\
\mathrm{transformation}.
\end{eqnarray}
Using $R_0^{-1}dx^1=dx_z^1+\om_0\times V_0^1$ and (\ref{eq:dify2})
we get $R_0^{-1}dx^1=dx_z^1-2(\pa_{v_1}x_z^1)^TN_0^0N_0^0dv_1$ and
\begin{eqnarray}\label{eq:linel}
|dx^1|^2=|dx_z^1|-4(\pa_{u_1}x_z^1)^TN_0^0(\pa_{v_1}x_z^1)^TN_0^0du^1\odot
dv^1,
\end{eqnarray}
so the linear element of the leaf $x^1$ does not depend on the
shape of the seed $x^0$ (we have ACLGN) provided
$(V_0^1)^TN_0^0=0\Rightarrow\mathcal{B}_1=\mathcal{B}_1(u_1,v_1)
(\Leftrightarrow d\mathcal{B}_1=0)$ (that is $\mathcal{B}_1$ does
not vary when $x_0^0$ varies on the tangent cone of $x_0$ from
$x_z$ for $x_z$ fixed). But for $x_z$ fixed and $x_0^0$ varying in
the TC $(V_0^1)^TN_0^0=0$ we have $(V_0^1)^TdN_0^0=0$ and
$d\log(\mathcal{B}_1)=z^{-1}\mathcal{B}_1(dN_0^0)^T(\pa_{u_1}x_z^1(\pa_{v_1}x_z^1)^T
+\pa_{v_1}x_z^1(\pa_{u_1}x_z^1)^T)N_0^0=^{(\ref{eq:int})}0$.

Thus all degenerate leaves $v_1=$ct ($u_1=$ct) are applicable to a
surface $y_0=y_0(u_1,v_1)$, the simple $\infty$ multiplicity of
facets of each leaf at a point $x_z$ corresponding to the tangent
planes of $y_0$ at a point counted with simple $\infty$
multiplicity and with $\pa_{u_1}x_z^1$ ($\pa_{v_1}x_z^1$) taken to
the same direction $\pa_{u_1}y_0$ ($\pa_{v_1}y_0$); by reflection
in $Tx_0$ (along the points of $x_0^0$ on the tangent cone of
$x_0$ from $x_z$) we get also the fact that
$(I_3-2N_0^0(N_0^0)^T)\pa_{v_1}x_z^1$
($(I_3-2N_0^0(N_0^0)^T)\pa_{u_1}x_z^1$) are taken to the same
direction, since the angle between these directions and the
initial ones $\pa_{u_1}x_z^1$ ($\pa_{v_1}x_z^1$) is independent of
$(u_0,v_0)$; thus these last directions must be taken to
$\pa_{v_1}y_0$ ($\pa_{u_1}y_0$).

For quadrics with $u_1,v_1=$ct being the ruling families on the
quadric $x_z$ confocal to $x_0$ (\ref{eq:int}) (which gives
integrability and applicability correspondence of a general
nature) is satisfied, but $m_0^1$ (respectively ${m_0^1}'$) are
defined by (\ref{eq:mm'}) with a-priori
$\mathcal{B}_1=\mathcal{B}_1(u_1,v_1)$ such that they depend only
on $x_0^0,v_1$ (respectively $x_0^0,u_1$) (they depend
quadratically on $v_1$ (respectively $u_1$), which makes
(\ref{eq:BB'}) a Ricatti equation). In this case by requiring that
$m_0^1,{m_0^1}'$ as defined by (\ref{eq:mm'}) satisfy the
reflection property (\ref{eq:mr}) we get the needed part of
(\ref{eq:int}) (without the coefficients of $V_0^1,N_0^0$) by
purely algebraic manipulations and the direct proof of the
integrability of (\ref{eq:BB'}) is simpler and enforces the
definition of $\mathcal{B}_1$ from (\ref{eq:mm'}): imposing the
compatibility condition $2z\ti d\wedge$ on the first equation of
(\ref{eq:BB'}) and using the equation itself we get the IC
$0=2z(\pa_{v_1}m_0^1dv_1+\mathcal{B}_1dx_0^0\times
\pa_{u_1}x_z^1)^T\wedge\om_0+2z(m_0^1)^Td\wedge\om_0=^{(\ref{eq:om})}
-(N_0^0)^T(2zm_0^1+m_0^1\times\pa_{v_1}m_0^1)(N_0^0)^T\frac{\om_0\times\wedge\om_0}{2}$;
using (\ref{eq:mr}), (\ref{eq:int}) and (\ref{eq:mm'}) this
becomes $0=(N_0^0)^T(2zm_0^1+m_0^1\times\pa_{v_1}m_0^1)=
\frac{z(m_0^1)^T\pa_{v_1}x_z^1}{(N_0^0)^T\pa_{v_1}x_z^1}-
\mathcal{B}_1(\pa_{u_1}x_z^1)^TN_0^0(V_0^1)^T\pa_{v_1}m_0^1=\\
\frac{z(V_0^1)^T(\mathcal{B}_1\pa_{v_1}x_z^1\times\pa_{u_1}x_z^1
+\pa_{v_1}(\mathcal{B}_1\pa_{u_1}x_z^1\times
V_0^1))}{(N_0^0)^T\pa_{v_1}x_z^1}$, which is straightforward;
replacing $(m_0^1,v_1)$ with $(m'^1_0,u_1)$ we get a similar
relation.

\subsection{The tangency configuration only symmetrizes the
integrability condition}\label{subsec:tcsym}

Consider first the case when the TC $(V_0^1)^TN_0^0=0$ is used
only to symmetrize (\ref{eq:int}):
\begin{eqnarray}\label{eq:int0}
M\pa_{u_1v_1}^2x_z-x_z=\mathrm{ct}(=0), \
(I_3-N_0N_0^T)[M(\pa_{u_1}x_z\pa_{v_1}x_z^T+\pa_{v_1}x_z\pa_{u_1}x_z^T)-x_zx_z^T
+x_0x_0^T]N_0=0\bigvee\nonumber\\
M\pa_{u_1v_1}^2x_z-x_z\neq\mathrm{ct},\
(I_3-N_0N_0^T)[\pa_{u_1}x_z\pa_{v_1}x_z^T+\pa_{v_1}x_z\pa_{u_1}x_z^T
-\pa_{u_1v_1}^2x_zx_z^T-x_z\pa_{u_1v_1}^2x_z^T+\nonumber\\
\pa_{u_1v_1}^2x_zx_0^T+x_0\pa_{u_1v_1}^2x_z^T]N_0=0.
\end{eqnarray}
In the first case with
$\mathcal{M}=\mathcal{M}(u_1,v_1):=M(\pa_{u_1}x_z\pa_{v_1}x_z^T+
\pa_{v_1}x_z\pa_{u_1}x_z^T)-x_zx_z^T$ we have
$(I_3-N_0N_0^T)[\mathcal{M}(u_1,v_1)-\mathcal{M}(\ti u_1,\ti
v_1)]N_0=0$ with $u_0,v_0,u_1,v_1,\ti u_1,\ti v_1$ independent
variables; applying $d$ we get
$[\mathcal{M}(u_1,v_1)-\mathcal{M}(\ti u_1,\ti v_1)]dN_0=
(I_3-N_0N_0^T)[\mathcal{M}(u_1,v_1)-\mathcal{M}(\ti u_1,\ti
v_1)]dN_0= N_0^T[\mathcal{M}(u_1,v_1)-\mathcal{M}(\ti u_1,\ti
v_1)]N_0dN_0$, so $\mathcal{M}(u_1,v_1)=\mathcal{M}(\ti u_1,\ti
v_1)+N_0^T[\mathcal{M}(u_1,v_1)-\mathcal{M}(\ti u_1,\ti
v_1)]N_0I_3$; differentiating this with respect to $u_1,v_1$ we
get $\pa_{u_1}(M\pa_{u_1}x_z)\pa_{v_1}x_z^T+
\pa_{v_1}x_z\pa_{u_1}(M\pa_{u_1}x_z)^T=N_0^T\pa_{u_1}\mathcal{M}(u_1,v_1)N_0I_3,\\
\pa_{v_1}(M\pa_{v_1}x_z)\pa_{u_1}x_z^T+
\pa_{u_1}x_z\pa_{v_1}(M\pa_{v_1}x_z)^T=N_0^T\pa_{v_1}\mathcal{M}(u_1,v_1)N_0I_3$,
so $\pa_{u_1}(M\pa_{u_1}x_z)=\pa_{v_1}(M\pa_{v_1}x_z)\\=0,\
\mathcal{M}=$ct and we get $x_z$ ((isotropic) singular) QC doubly
ruled by leaves (the image of the unit sphere
$X=X(u_1,v_1):=\frac{-u_1v_1f_1+2\bar f_1+(u_1+v_1)e_3}{u_1-v_1}$
under a linear transformation of $\mathbb{C}^3$ with at most
$1$-dimensional kernel).

In the second case if $M\pa_{u_1v_1}^2x_z\neq w,\
w=\mathrm{ct}\in\mathbb{C}^3\setminus\{0\}$, then we shall see
later that the condition imposed on $x_0$ is over-determined.

Thus $M\pa_{u_1v_1}^2x_z=w,\
w=\mathrm{ct}\in\mathbb{C}^3\setminus\{0\}$; with
$\mathcal{M}:=M(\pa_{u_1}x_z\pa_{v_1}x_z^T+\pa_{v_1}x_z\pa_{u_1}x_z^T)
-x_zw^T-wx_z^T$ and as above we get
$\pa_{u_1}(M\pa_{u_1}x_z)=\pa_{v_1}(M\pa_{v_1}x_z)=0,\
\mathcal{M}=$ct and we get $x_z$ ((isotropic) singular) (I)QWC
doubly ruled by leaves (the image of the equilateral paraboloid
$Z=Z(u_1,v_1):=u_1f_1+v_1\bar f_1+u_1v_1e_3$ under an affine
(which can be taken to be linear) transformation of $\mathbb{C}^3$
with kernel $\neq\mathbb{C}e_3$ at most $1$-dimensional).

With $J_2:=f_1f_1^T,\ J_3:=f_1e_3^T+e_3f_1^T$ being the {\it
symmetric Jordan} (SJ) blocks of dimension $2$ and $3$ any
$A=A^T\in\mathbf{M}_3(\mathbb{C})$ can be brought, via conjugation
with a complex rotation, to a form with diagonal blocks, each
block being a number or $aI_2+J_2$ or $aI_3+J_3,\ a\in\mathbb{C}$
and, since one can take square roots of any SJ matrix without
isotropic kernel, any $A\in\mathbf{GL}_3(\mathbb{C})$ has {\it
singular value decomposition} (SVD) $A=R_1DR_2^T,\
R_1,R_2\in\mathbf{O}_3(\mathbb{C}),\ D^2$ SJ.

We have
\begin{eqnarray}\label{eq:XX}
-\frac{1}{2}(u-v)^2(\pa_uX\pa_vX^T+\pa_vX\pa_uX^T)-XX^T=-I_3,\
\pa_{uv}^2X=-\frac{2}{(u-v)^2}X,\ X=X(u,v)\nonumber\\
\pa_uZ\pa_vZ^T+\pa_vZ\pa_vZ^T
-\pa_{uv}^2ZZ^T-Z\pa_{uv}^2Z^T=f_1\bar f_1^T+\bar f_1f_1^T,\
Z=Z(u,v),
\end{eqnarray}
so (\ref{eq:int0}) becomes
\begin{eqnarray}\label{eq:int1}
(I_3-N_0N_0^T)(-\mathcal{A}\mathcal{A}^T+x_0x_0^T)N_0=0,\
\dim(\ker(\mathcal{A}))\le 1\
\mathrm{for}\ x_z=\mathcal{A}X_1 \mathrm{((isotropic)\ singular)\ QC},\nonumber\\
(I_3-N_0N_0^T)(\mathcal{A}(f_1\bar f_1^T+\bar
f_1f_1^T)\mathcal{A}^T
+x_0(\mathcal{A}e_3)^T+\mathcal{A}e_3x_0^T)N_0=0,\
\dim(\ker(\mathcal{A}))\le 1,\nonumber\\
\mathcal{A}e_3\neq 0\ \mathrm{for}\ x_z=\mathcal{A}Z_1\
\mathrm{((isotropic)\ singular)\ (I)QWC};
\end{eqnarray}
for $x_z$ ((isotropic) singular) QC one can multiply $\mathcal{A}$
on the right with a rotation $R_2\in\mathbf{O}_3(\mathbb{C})$
(which has the effect of changing each of $u_1,v_1$ by the same
M\"{o}bius transformation); for $x_z$ ((isotropic) singular)
(I)QWC one can multiply $\mathcal{A}$ on the right with
$e^{c_1}(e^{ic_2}f_1\bar f_1^T+e^{-ic_2}\bar
f_1f_1^T)+e^{2c_1}e_3e_3^T$ (which has the effect
$(u_1,v_1)\rightarrow(e^{-c_1-ic_2}u_1,e^{-c_1+ic_2}v_1)$) or with
the reflection $e_2\mapsto -e_2$ (which has the effect
$u_1\leftrightarrow v_1$); in both cases one can multiply
$\mathcal{A}$ on the left with an arbitrary rotation (which must
also apply to $x_0$).

For $x_z$ ((isotropic) singular) QC apply this rotation to bring
$\mathcal{A}\mathcal{A}^T$ to the SJ canonical form; for
$\ker(\mathcal{A}\mathcal{A}^T)$ containing the isotropic
direction $f_1$ we have $\mathcal{A}\mathcal{A}^T=J_2+ae_3e_3^T,\
a\in\mathbb{C}$ or $=J_3$ or $=ae_3e_3^T,\
a\in\mathbb{C}\setminus\{0\}$.

In the first case with $x_0=:\frac{\mathbf{w}}{\sqrt{2}}f_1+
\frac{\mathbf{\ti w}}{\sqrt{2}}\bar f_1+\sqrt{2\mathbf{z}}e_3$ we
have
$dx_0=(\frac{1}{\sqrt{2}}f_1+\frac{\pa_{\mathbf{w}}\mathbf{z}}
{\sqrt{2\mathbf{z}}}e_3)d\mathbf{w}+(\frac{1}{\sqrt{2}}\bar
f_1+\frac{\pa_{\mathbf{\ti w}}\mathbf{z}}
{\sqrt{2\mathbf{z}}}e_3)d\mathbf{\ti w},\
N_0=\frac{\pa_{\mathbf{\ti
w}}\mathbf{z}f_1+\pa_{\mathbf{w}}\mathbf{z}\bar
f_1-\frac{\sqrt{2\mathbf{z}}}{\sqrt{2}}e_3}
{\sqrt{2\pa_{\mathbf{w}}\mathbf{z}\pa_{\mathbf{\ti
w}}\mathbf{z}+\mathbf{z}}}$ and we need
$2a\pa_{\mathbf{w}}\mathbf{z} +(\mathbf{\ti
w}+2\pa_{\mathbf{w}}\mathbf{z})(\mathbf{w}\pa_{\mathbf{w}}\mathbf{z}+
\mathbf{\ti w}\pa_{\mathbf{\ti
w}}\mathbf{z}-2\mathbf{z})=-2\pa_{\mathbf{w}}\mathbf{z}+2a\pa_{\mathbf{\ti
w}}\mathbf{z} +(\mathbf{w}+2\pa_{\mathbf{\ti
w}}\mathbf{z})(\mathbf{w}\pa_{\mathbf{w}}\mathbf{z}+ \mathbf{\ti
w}\pa_{\mathbf{\ti w}}\mathbf{z}-2\mathbf{z})=0,\
2\pa_{\mathbf{w}}\mathbf{z}\pa_{\mathbf{\ti
w}}\mathbf{z}+\mathbf{z}\neq 0$.

For $a\neq 0$ we have $\mathbf{\ti
w}+2\pa_{\mathbf{w}}\mathbf{z}\neq 0$, so $\pa_{\mathbf{\ti
w}}\mathbf{z}=-\frac{2a\pa_{\mathbf{w}}\mathbf{z} +(\mathbf{\ti
w}+2\pa_{\mathbf{w}}\mathbf{z})(\mathbf{w}\pa_{\mathbf{w}}\mathbf{z}
-2\mathbf{z})}{\mathbf{\ti w}(\mathbf{\ti
w}+2\pa_{\mathbf{w}}\mathbf{z})}$ and $\pa_{\mathbf{w}}\mathbf{z}$
is a root of the cubic equation $\mathcal{X}^3+(\mathbf{\ti
w}+a\mathbf{w})\mathcal{X}^2+\frac{1}{4}(\mathbf{\ti
w}^2+2a\mathbf{w}\mathbf{\ti
w}-4a\mathbf{z}+2a^2)\mathcal{X}-\frac{1}{2}a\mathbf{\ti
w}\mathbf{z}=0$ in $\mathcal{X}$; the condition that this equation
has a triple root is over-determined (it imposes two functionally
independent conditions on $\mathbf{w},\mathbf{\ti w},\mathbf{z}$)
and the condition of it having two distinct roots (with
$\pa_{\mathbf{w}}\mathbf{z}$ double root) leads to $x_0$ being the
isotropic developable circumscribed to one (and thus to the other)
of the conics
$\frac{\mathbf{w}^2}{2}+\frac{1}{a}2\mathbf{z}=1\bigwedge\mathbf{\ti
w}=0,\ -\frac{1}{2a^2}\mathbf{\ti
w}^2-\frac{1}{a}\mathbf{w}\mathbf{\ti w}=1\bigwedge\mathbf{z}=0$
which extended with the isotropic plane $f_1^Tx=0$ determines the
confocal quadrics below.

Otherwise applying $\pa_{\mathbf{w}}$ to this cubic equation we
get $\pa_{\mathbf{w}}^2\mathbf{z}=0$, so
$\mathbf{z}=\mathbf{w}F(\mathbf{\ti w})+G(\mathbf{\ti w})$ and we
need $0=[\mathbf{\ti w}+2F(\mathbf{\ti w})][\mathbf{\ti
w}F'(\mathbf{\ti w})-F(\mathbf{\ti w})]=2aF(\mathbf{\ti
w})+[\mathbf{\ti w}+2F(\mathbf{\ti w})][\mathbf{\ti
w}G'(\mathbf{\ti w})-2G(\mathbf{\ti w})]=[1+2F'(\mathbf{\ti
w})][\mathbf{\ti w}F'(\mathbf{\ti w})-F(\mathbf{\ti w})]=
2aF'(\mathbf{\ti w})+[1+2F'(\mathbf{\ti w})][\mathbf{\ti
w}G'(\mathbf{\ti w})-2G(\mathbf{\ti w})]+G'(\mathbf{\ti
w})[\mathbf{\ti w}F'(\mathbf{\ti w})-F(\mathbf{\ti
w})]=-2F(\mathbf{\ti w})+G'(\mathbf{\ti w})[2a+\mathbf{\ti
w}G'(\mathbf{\ti w})-2G(\mathbf{\ti w})]$, so $F(\mathbf{\ti
w})=c\mathbf{\ti w},\
c\in\mathbb{C}\setminus\{-1,-\frac{1}{2},0\},\ G(\mathbf{\ti
w})=\frac{ac}{1+2c}+\frac{c(1+2c)}{2a(1+c)}\mathbf{\ti w}^2$; this
leads to $x_0$ being one of the confocal Darboux quadrics
$x_w^TAR_w^{-1}x_w=1,\ A:=(a_1I_3+J_2+ae_3e_3^T)^{-1},\
a_1\in\mathbb{C}\setminus\{0,-a\},\ R_w:=I_3-wA,\
w\in\mathbb{C}\setminus\{a_1,a_2:=a_1+a\}$; since
$AR_w^{-1}=((a_1-w)I_3+J_2+ae_3e_3^T)^{-1}$ we can take $w=0$ for
$x_0$. In this case
$\mathcal{A}_z=f_1(\frac{z}{2}f_1+\frac{1}{z}\bar
f_1)^T+\sqrt{a}e_3e_3^T,\ z\in\mathbb{C}\setminus\{0\}$ and
$x_z=\frac{z^2-u_1v_1}{z(u_1-v_1)}f_1+\sqrt{a}\frac{u_1+v_1}{u_1-v_1}e_3,\
z\in\mathbb{C}\setminus\{0\}$ is the isotropic plane $f_1^Tx_z=0$
with rulings the lines tangent to the conic $(\bar
f_1^Tx_z)^2+\frac{1}{a}(e_3^Tx_z)^2=1$ obtained for $z^2+u_1v_1=0$
(this is the same conic as the one in the plane $\mathbf{\ti w}=0$
generating the isotropic developable above); this is the isotropic
singular quadric of the above considered family of confocal
quadrics (obtained for $w=a_1$).

For $x_0^0\subset x_0$ from (\ref{eq:XX}) and
$(x_z^1)^TN_0^0=(x_0^0)^TN_0^0$ we get
$4(\pa_{u_1}x_z^1)^TN_0^0(\pa_{v_1}x_z^1)^TN_0^0du_1\odot
dv_1=-\frac{4a_1du_1\odot dv_1}{(u_1-v_1)^2}=-a_1|dX_1|^2$ and
from (\ref{eq:linel}) the linear element of the leaves is
$|dx^1|^2=|dx_z^1|^2+a_1|dX_1|^2=
dX_1^T(\mathcal{A}^T\mathcal{A}+a_1I_3)dX_1=dX_1^T(a_1I_3+ae_3e_3^T)dX_1$,
so the leaves are applicable to a quadric with center and of
revolution.

For $a=0$ $\mathbf{w}\pa_{\mathbf{w}}\mathbf{z}+ \mathbf{\ti
w}\pa_{\mathbf{\ti w}}\mathbf{z}-2\mathbf{z}=0$ leads to $x_0$
being the plane passing through the origin and having unit normal
$N_0=cf_1+e_3,\ c\in\mathbb{C}$ and $\mathbf{\ti
w}+2\pa_{\mathbf{w}}\mathbf{z}=0$ leads to
$\mathbf{z}+\frac{\mathbf{w}\mathbf{\ti w}}{2}=F(\mathbf{\ti w}),\
F'(\mathbf{\ti w})^2-2\frac{F(\mathbf{\ti w})}{\mathbf{\ti
w}}F'(\mathbf{\ti w})+\frac{1}{2}=0$ with solutions $F(\mathbf{\ti
w})=\pm\frac{1}{\sqrt{2}}\mathbf{\ti w}$, which leads to $x_0$
being the isotropic cones centered at $\pm f_1$ (the isotropic
developable which extended with the isotropic plane $f_1^Tx=0$
generates the confocal quadrics below) and $F(\mathbf{\ti
w})=c\mathbf{\ti w}^2+\frac{1}{8c},\
c\in\mathbb{C}\setminus\{0\}$, which leads to $x_0$ being one of
the confocal Darboux quadrics $x_w^TAR_w^{-1}x_w=1,\
A:=(a_1I_3+J_2)^{-1},\ a_1\in\mathbb{C}\setminus\{0\},\
w\in\mathbb{C}\setminus\{a_1\}$ (the only quadrics having contact
of order $3$ with $C(\infty)$); we can take $w=0$ for $x_0$. In
this case $\mathcal{A}_z=f_1e_3^T+(-\frac{z^2}{2}f_1+\bar
f_1+ze_3)f_1^T,\
x_z=\frac{u_1+v_1}{u_1-v_1}f_1+\frac{2}{u_1-v_1}(-\frac{z^2}{2}f_1+\bar
f_1+ze_3),\ z\in\mathbb{C},\ \mathcal{A}_{\infty}=J_3,\
x_{\infty}=\frac{u_1+v_1}{u_1-v_1}f_1+\frac{2}{u_1-v_1}e_3$ (after
multiplication with a rotation $R_2\in\mathbf{O}_3(\mathbb{C})$ on
the right one can make $\mathcal{A}_z=f_1e_3^T+(\bar
f_1+ze_3)f_1^T$ and thus the definition of $\mathcal{A}_{\infty}$
is clear) is the pencil of (isotropic) planes containing $f_1$ and
the rulings are lines passing through $\pm f_1$ in their
respective planes (these foci are obtained from
$x_z(u_1,\infty),x_z(\infty,v_1)$); the only finite (isotropic)
singular quadric of the confocal family of $x_0$ is the line
$x_{a_1}=\mathbb{C}f_1$ which is also singular set of $x_z,\
z\in\mathbb{C}\cup\{\infty\}$; however for convenience one can
also include the pencil of (isotropic) planes containing $f_1$ in
this isotropic singular quadric.

From (\ref{eq:XX}) and $(x_z^1)^TN_0^0=(x_0^0)^TN_0^0$ we get
$4(\pa_{u_1}x_z^1)^TN_0^0(\pa_{v_1}x_z^1)^TN_0^0du_1\odot
dv_1=-a_1|dX_1|^2$ and from (\ref{eq:linel}) the linear element of
the leaves is $|dx^1|^2=|dx_z^1|^2+a_1|dX_1|^2=
dX_1^T(\mathcal{A}_z^T\mathcal{A}_z+a_1I_3)dX_1$; since
$\mathcal{A}_z^T\mathcal{A}_z=J_3$ for $z\in\mathbb{C},\
=J_3^2=J_2$ for $z=\infty$, the leaves are applicable to another
region $x_0^1$ of $x_0$ for $z=\infty$ and for $z\in\mathbb{C}$
they are applicable to a Darboux quadric with $A^{-1}=a_1I_3+J_3$
(thus different from the type of $x_0$). In this case
(\ref{eq:BB'}) are linear ($m_0^1,{m_0^1}'$ depend linearly on
$v_1,u_1$), so this B transformation can be found by quadratures.

For $\mathcal{A}\mathcal{A}^T=J_3$ with
$x_0=:\frac{\mathbf{w}}{\sqrt{2}}f_1+ \frac{\mathbf{\ti
w}}{\sqrt{2}}\bar f_1+\mathbf{z}e_3$ we have
$dx_0=\frac{1}{\sqrt{2}}(\pa_{\mathbf{\ti w}}\mathbf{w}f_1+\bar
f_1)d\mathbf{\ti
w}+(\frac{\pa_{\mathbf{z}}\mathbf{w}}{\sqrt{2}}f_1+e_3)d\mathbf{z},\
N_0=\frac{\sqrt{2}(\pa_{\mathbf{\ti w}}\mathbf{w}f_1-\bar
f_1)+\pa_{\mathbf{z}}{\mathbf{w}}e_3}
{\sqrt{(\pa_{\mathbf{z}}\mathbf{w})^2-4\pa_{\mathbf{\ti
w}}\mathbf{w}}}$ and we need $-\sqrt{2}\pa_{\mathbf{z}}\mathbf{w}
+(\mathbf{\ti w}\pa_{\mathbf{\ti
w}}\mathbf{w}+\mathbf{w})(-\mathbf{w}+\mathbf{\ti
w}\pa_{\mathbf{\ti
w}}\mathbf{w}+\mathbf{z}\pa_{\mathbf{z}}\mathbf{w})=
2\sqrt{2}+(2\mathbf{z}+\mathbf{\ti
w}\pa_{\mathbf{z}}\mathbf{w})(-\mathbf{w}+\mathbf{\ti
w}\pa_{\mathbf{\ti
w}}\mathbf{w}+\mathbf{z}\pa_{\mathbf{z}}\mathbf{w})=0,\
(\pa_{\mathbf{z}}\mathbf{w})^2-4\pa_{\mathbf{\ti w}}\mathbf{w}\neq
0$. We have $2\mathbf{z}+\mathbf{\ti
w}\pa_{\mathbf{z}}\mathbf{w}\neq 0,\ \pa_{\mathbf{\ti
w}}\mathbf{w}=-\frac{2\sqrt{2}+(2\mathbf{z}+\mathbf{\ti
w}\pa_{\mathbf{z}}\mathbf{w})(-\mathbf{w}
+\mathbf{z}\pa_{\mathbf{z}}\mathbf{w})}{\mathbf{\ti
w}(2\mathbf{z}+\mathbf{\ti w}\pa_{\mathbf{z}}\mathbf{w})}$ and
$\pa_{\mathbf{z}}\mathbf{w}$ solution of the cubic equation
$\mathcal{X}^3+\frac{2\mathbf{z}}{\mathbf{\ti
w}}\mathcal{X}^2+\frac{4\mathbf{w}}{\mathbf{\ti
w}}\mathcal{X}+\frac{4(2\mathbf{w}\mathbf{z}-\sqrt{2})}{\mathbf{\ti
w}^2}=0$ in $\mathcal{X}$; the condition that this cubic equation
has a triple root is over-determined and the condition that of it
having two distinct roots (with $\pa_{\mathbf{z}}\mathbf{w}$
double root) leads to $x_0$ being the isotropic developable
circumscribed to the conic
$\sqrt{2}\mathbf{w}\mathbf{z}=1\bigwedge\mathbf{\ti w}=0$ which
extended with the isotropic plane $f_1^Tx=0$ determines the
confocal quadrics below.

Otherwise applying $\pa_{\mathbf{z}}$ to this equation we get
$\pa_{\mathbf{z}}^2\mathbf{w}=-\frac{2}{\mathbf{\ti w}}$, so
$\mathbf{w}=-\frac{\mathbf{z}^2}{\mathbf{\ti
w}}+\mathbf{z}F(\mathbf{\ti w})+G(\mathbf{\ti w})$ and we need
$0=F(\mathbf{\ti w})F'(\mathbf{\ti w})=F'(\mathbf{\ti w})^2$, so
$F(\mathbf{\ti w})=c,\ c\in\mathbb{C}$ and further
$0=2\sqrt{2}+c\mathbf{\ti w}[\mathbf{\ti w}G'(\mathbf{\ti
w})-G(\mathbf{\ti w})]=-\sqrt{2}c+[\mathbf{\ti w}G'(\mathbf{\ti
w})+G(\mathbf{\ti w})][\mathbf{\ti w}G'(\mathbf{\ti
w})-G(\mathbf{\ti w})]$, so $c\neq 0$ and $G(\mathbf{\ti
w})=-\frac{c^2}{2}\mathbf{\ti w}+\frac{2\sqrt{2}}{c\mathbf{\ti
w}}$; this leads to $x_0$ being one of the confocal Darboux
quadrics $x_w^TAR_w^{-1}x_w=1,\ A:=(a_1I_3+J_3)^{-1},\
a_1\in\mathbb{C}\setminus\{0\},\ w\in\mathbb{C}\setminus\{a_1\}$;
we can take $w=0$ for $x_0$. In this case
$\mathcal{A}_z=f_1(-\frac{z^2}{2}f_1+\bar f_1+ze_3)^T+e_3f_1^T,\
x_z=-\frac{(u_1-z)(v_1-z)}{u_1-v_1}f_1+\frac{2}{u_1-v_1}e_3,\
z\in\mathbb{C}$ is the isotropic plane $f_1^Tx_z=0$ with rulings
the lines tangent to the conic $2\bar f_1^Tx_ze_3^Tx_z=1$ obtained
for $(u_1-z)+(v_1-z)=0$; this is the isotropic singular quadric of
the above considered family of confocal quadrics (obtained for
$w=a_1$); as above the leaves are applicable to the Darboux
quadric with $A^{-1}=a_1I_3+J_2$ discussed above, so this B
transformation is the inversion of the previous one.

For $\mathcal{A}\mathcal{A}^T=ae_3e_3^T,\
a\in\mathbb{C}\setminus\{0\}$ with
$x_0=:\sqrt{\mathbf{x}}e_1+\sqrt{\mathbf{y}}e_2+\sqrt{\mathbf{z}}e_3$
we need
$\mathbf{x}\pa_{\mathbf{x}}\mathbf{z}+\mathbf{y}\pa_{\mathbf{y}}\mathbf{z}-\mathbf{z}
=-a\frac{\pa_{\mathbf{x}}\mathbf{z}}{1+\pa_{\mathbf{x}}\mathbf{z}}=
-a\frac{\pa_{\mathbf{y}}\mathbf{z}}{1+\pa_{\mathbf{y}}\mathbf{z}}$,
so $\pa_{\mathbf{x}}\mathbf{z}=\pa_{\mathbf{y}}\mathbf{z},\
\mathbf{z}=\mathbf{z}(\mathbf{x}+\mathbf{y})$ and
$\pa_{\mathbf{x}}\mathbf{z}$ is a root of the quadratic equation
$(1+\mathcal{X})[(\mathbf{x}+\mathbf{y})\mathcal{X}-\mathbf{z}
+a\frac{\mathcal{X}}{1+\mathcal{X}}]=0$ in $\mathcal{X}$; the
equation has a double root for the isotropic developable
$(\mathbf{x}+\mathbf{y}+\mathbf{z}-a)^2=-4a(\mathbf{x}+\mathbf{y})\
(\Leftrightarrow\
\mathbf{x}+\mathbf{y}+(\sqrt{\mathbf{z}}\pm\sqrt{a})^2=0$, that is
the isotropic cones centered at $\pm\sqrt{a}e_3$) circumscribed to
the circle $\mathbf{x}+\mathbf{y}+a=0\bigwedge\mathbf{z}=0$, which
extended with the isotropic planes $\mathbf{x}+\mathbf{y}=0$ (for
which the above quadratic equation becomes linear) determines the
confocal quadrics of revolution below.

Otherwise applying $\pa_\mathbf{x}$ to this equation we get
$\pa_{\mathbf{x}}^2\mathbf{z}=0$, so
$\mathbf{z}(\mathbf{x}+\mathbf{y})=c(1+\frac{\mathbf{x}+\mathbf{y}}{a-c}),\
c\in\mathbb{C}\setminus\{0,a\}$, which gives the family of
confocal quadrics of revolution $x_w^TAR_{w}^{-1}x_w=1,\
A:=(a_1I_3+ae_3e_3^T)^{-1},\
a_1\in\mathbb{C}\setminus\{0,-a\},w\in\mathbb{C}\setminus\{a_1,\
a_2:=a_1+a\}$. In this case
$\mathcal{A}_z=(\frac{z}{2}f_1+\frac{1}{z}\bar
f_1)f_1^T+\sqrt{a}e_3e_3^T,\
x_z=\frac{1}{u_1-v_1}(zf_1+\frac{2}{z}\bar
f_1)+\sqrt{a}\frac{u_1+v_1}{u_1-v_1}e_3,\
z\in\mathbb{C}\setminus\{0\},\ \mathcal{A}_0=\bar
f_1f_1^T+\sqrt{a}e_3e_3^T,\ x_{z=0}=\frac{2}{u_1-v_1}\bar
f_1+\sqrt{a}\frac{u_1+v_1}{u_1-v_1}e_3,\
\mathcal{A}_{\infty}=f_1f_1^T+\sqrt{a}e_3e_3^T,\
x_{\infty}=\frac{2}{u_1-v_1}f_1+\sqrt{a}\frac{u_1+v_1}{u_1-v_1}e_3$
is the pencil of (isotropic) planes containing $e_3$ and the
rulings are lines passing through $\pm\sqrt{a}e_3$ in their
respective planes (these foci are obtained from
$x_z(u_1,\infty),x_z(\infty,v_1)$); in this case the isotropic
singular quadric of the confocal family is the isotropic planes
$f_1^Tx_{a_1}=0,\bar f_1^Tx_{a_1}=0$; its singular part is the
$e_3$ axis which appears as the singular set of $x_z,\
z\in\mathbb{C}\cup\{\infty\}$ and gives the complementary
transformation; however for convenience one can also include the
pencil of planes containing $e_3$ in this isotropic singular
quadric (they form a triply orthogonal system with the remaining
quadrics of the confocal family).

We have $\mathcal{A}_z^T\mathcal{A}_z=J_2+ae_3e_3^T,\
z\in\mathbb{C}\setminus\{0\},\
\mathcal{A}_0^T\mathcal{A}_0=\mathcal{A}_{\infty}^T\mathcal{A}_{\infty}=ae_3e_3^T$;
thus as above for $z=0,\infty$ the leaves are applicable to
another region $x_0^1$ of $x_0$ and for
$z\in\mathbb{C}\setminus\{0\}$ the leaves are applicable to the
Darboux quadric with $A^{-1}=a_1I_3+J_2+ae_3e_3^T,\
a_1\in\mathbb{C}\setminus\{0,-a\}$ discussed above; thus this B
transformation is the inversion of the previous one.

For $\ker(\mathcal{A}\mathcal{A}^T)=\mathbb{C}e_3$ we have
$\mathcal{A}\mathcal{A}^T=a(f_1\bar f_1^T+\bar f_1f_1^T)+J_2,\
a\in\mathbb{C}\setminus\{0\}$ or $=ae_1e_1^T+be_2e_2^T,\
a,b\in\mathbb{C}\setminus\{0\}$.

Since by adding a multiple of $I_3$ to $\mathcal{A}\mathcal{A}^T$
does not change (\ref{eq:int1}), the first case and the second one
for $a=b$ have already been discussed for $x_0$; the only change
is for $x_z$: we obtain the singular $B_z$ transformation, when
the applicability correspondence is given by the Ivory affinity
between confocal quadrics
$(\sqrt{A})^{-1}X\rightarrow\sqrt{R_z}(\sqrt{A})^{-1}X$ (which can
be extended to an affine correspondence between $x_0$ and the
singular $x_z$ since one can take square roots of symmetric
matrices with non-isotropic kernels).

For $a\neq b$ with
$x_0=:\sqrt{\mathbf{x}}e_1+\sqrt{\mathbf{y}}e_2+\sqrt{\mathbf{z}}e_3$
we need $\pa_\mathbf{y}\mathbf{z}=
\frac{a\pa_\mathbf{x}\mathbf{z}}{b-(a-b)\pa_\mathbf{x}\mathbf{z}}$
and $\pa_\mathbf{x}\mathbf{z}$ solution of the cubic equation
$(1+\mathcal{X})(b-(a-b)\mathcal{X})[\mathbf{x}\mathcal{X}+
\mathbf{y}\frac{a\mathcal{X}}{b-(a-b)\mathcal{X}}-\mathbf{z}-
a\frac{\mathcal{X}}{1+\mathcal{X}}]=0$ in $\mathcal{X}$; the
condition that this cubic equation has a triple root is
over-determined and the condition of it having two distinct roots
(with $\pa_\mathbf{x}\mathbf{z}$ double root) leads to $x_0$ being
the isotropic developable circumscribed to one (and thus to all)
of the conics
$\frac{\mathbf{x}}{a}+\frac{\mathbf{y}}{b}=1\bigwedge\mathbf{z}=0,\
\frac{\mathbf{y}}{b-a}-\frac{\mathbf{z}}{a}=1\bigwedge\mathbf{x}=0,\
\frac{\mathbf{x}}{a-b}-\frac{\mathbf{z}}{b}=1\bigwedge\mathbf{y}=0$
which determines the confocal quadrics below.

Otherwise applying $\pa_\mathbf{y}$ to the cubic equation and
replacing $\pa_\mathbf{y}\mathbf{z}$ with its value we get
$\pa_\mathbf{y}(\pa_\mathbf{x}\mathbf{z})=0$, so $\mathbf{z}$ is
linear in $\mathbf{x},\mathbf{y}$, which gives the family of
confocal quadrics $x_w^TAR_{w}^{-1}x_w=1,\
A^{-1}:=(a+a_3)e_1e_1^T+(b+a_3)e_2e_2^T+a_3e_3e_3^T,\
a_3\in\mathbb{C}\setminus\{0,-a,-b\},\
w\in\mathbb{C}\setminus\{a_1:=a+a_3,a_2:=b+a_3,a_3\}$; we can take
$w=0$ for $x_0$. All confocal quadrics cut the above isotropic
developable along $4$ isotropic rulings of each ruling family and
passing through $(\frac{\mathbf{x}}{a}+\frac{\mathbf{y}}{b}=1,\
\mathbf{z}=0)\bigcap(\frac{\mathbf{x}}{a+a_3}+\frac{\mathbf{y}}{b+b_3}=1,\
\mathbf{z}=0)$ (any $3$ rulings of a family individuate the
quadric) which intersect in the $12$ finite umbilics and in $4$
points situated on $C(\infty)$ (Darboux).

We have $\mathcal{A}=\sqrt{a}e_1e_1^T+\sqrt{b}e_2e_2^T,\ x_z$ is
the singular quadric $x_{a_3}$ of the confocal family (the plane
$e_3^Tx=0$) and we have the singular $B_{a_3}$ transformation.

For $\ker(\mathcal{A}\mathcal{A}^T)=0$ all cases of quadrics
except the (pseudo-)sphere have been covered for $x_0$; using the
SVD of $\mathcal{A}$ we get $x_z$ quadric confocal to $x_0$ and
the usual $B_z$ transformation; the case of the (pseudo-)sphere
follows from $dx_0^Tx_0=0$.

For $x_z$ ((isotropic) singular) (I)QWC with $\mathcal{A}(f_1\bar
f_1^T+\bar f_1f_1^T)\mathcal{A}^T=0$ we have $\mathcal{A}f_1=0$;
multiplying $\mathcal{A}$ on the left with a rotation we get
$\mathcal{A}e_3=ae_3,\ a\in\mathbb{C}\setminus\{0\}$ or $=f_1$.

In the first case with
$x_0=:\sqrt{\mathbf{x}}e_1+\sqrt{\mathbf{y}}e_2+\mathbf{z}e_3$ we
need
$\mathbf{x}\pa_\mathbf{x}\mathbf{z}+\mathbf{y}\pa_\mathbf{y}\mathbf{z}-\mathbf{z}=
\frac{1}{4\pa_\mathbf{x}\mathbf{z}}=\frac{1}{4\pa_\mathbf{y}\mathbf{z}}$,
so $\mathbf{z}=\mathbf{z}(\mathbf{x}+\mathbf{y})$ and
$\pa_\mathbf{x}\mathbf{z}$ is solution of the quadratic equation
$4(\mathbf{x}+\mathbf{y})\mathcal{X}^2-4\mathbf{z}\mathcal{X}-1=0$
in $\mathcal{X}$; the equation has a double root for the isotropic
cone $\mathbf{x}+\mathbf{y}+\mathbf{z}^2=0$ which extended with
the isotropic planes $\mathbf{x}+\mathbf{y}=0$ (for which the
above quadratic equation becomes linear) determines the confocal
quadrics below.

Otherwise applying $\pa_\mathbf{x}$ to this equation we get
$\pa_\mathbf{x}^2\mathbf{z}=0$, so $x_0$ is one of the confocal
paraboloids of revolution
$\frac{\mathbf{x}+\mathbf{y}}{a_1-w}=2\mathbf{z}+(a_1-w),\
a_1\in\mathbb{C}\setminus\{0\},\ w\in\mathbb{C}\setminus\{a_1\}$;
we can take $w=0$ for $x_0$ (note that to obtain the canonical
form for $x_0$ we need to further apply the translation
$\frac{a_1}{2}e_3$). In this case we have
$\mathcal{A}=\mathbf{v}f_1^T+ae_3e_3^T,\ \mathbf{v}\times e_3\neq
0,\ x_z=v_1(\mathbf{v}+au_1e_3)$; by a transformation of $u_1,v_1$
into linear functions of themselves we can make
$\mathcal{A}_z=(\frac{z}{2}f_1+\frac{1}{z}\bar
f_1)f_1^T+e_3e_3^T,\ x_z=v_1(zf_1+\frac{2}{z}\bar f_1+u_1e_3),\
z\in\mathbb{C}\setminus\{0\},\ \mathcal{A}_0=\bar
f_1f_1^T+e_3e_3^T,\ x_{z=0}=v_1(\bar f_1+u_1e_3),\
\mathcal{A}_{\infty}=f_1f_1^T+e_3e_3^T,\
x_{\infty}=v_1(f_1+u_1e_3)$ the pencil of (isotropic) planes
containing $e_3$ and the rulings are lines passing through $0$ and
$\infty e_3$ in their respective planes; in this case the
isotropic singular quadric of the confocal family is the isotropic
planes $f_1^Tx_{a_1}=0,\bar f_1^Tx_{a_1}=0$; its singular part is
the $e_3$ axis which appears as the singular set of $x_z,\
z\in\mathbb{C}\cup\{\infty\}$ and gives the complementary
transformation; however for convenience one can also include the
pencil of planes containing $e_3$ in this isotropic singular
quadric (they form a triply orthogonal system with the remaining
quadrics of the confocal family).

For $x_0^0\subset x_0$ from (\ref{eq:XX}) and
$(x_z^1)^TN_0^0=(x_0^0)^TN_0^0$ we get
$4(\pa_{u_1}x_z^1)^TN_0^0(\pa_{v_1}x_z^1)^TN_0^0du_1\odot
dv_1=-2a_1du_1\odot dv_1=-a_1|(f_1\bar f_1^T+\bar
f_1f_1^T)dZ_1|^2$ and from (\ref{eq:linel}) the linear element of
the leaves is $|dx^1|^2=|dx_z^1|^2+a_1|(f_1\bar f_1^T+\bar
f_1f_1^T)dZ_1|^2= dZ_1^T[\mathcal{A}_z^T\mathcal{A}_z+a_1(f_1\bar
f_1^T+\bar f_1f_1^T)]dZ_1$; since
$\mathcal{A}_z^T\mathcal{A}_z=J_2+e_3e_3^T$ for
$z\in\mathbb{C}\setminus\{0\}$ and $=e_3e_3^T$ otherwise the
leaves are applicable to a Darboux quadric without center for
$z\in\mathbb{C}\setminus\{0\}$ and to another region of $x_0$
otherwise.

The case $\mathcal{A}e_3=f_1$ leads to $x_0$ being the plane
passing through the origin and having unit normal $N_0=cf_1+e_3,\
c\in\mathbb{C}$ or the isotropic cone $|x_0|^2=0$.

The remaining cases of IQWC should follow by similar computations.

\subsection{Rigidity of the B\"{a}cklund transformation of quadrics}

Given the defining surface $x_0$ being a quadric, we should get
the auxiliary surface $x_z$ as in \S \ref{subsec:tcsym}.

For $x_0=(\sqrt{A})^{-1}X(u_0,v_0)$ canonical QC with
$Y(v):=-v^2f_1+2\bar f_1+2ve_3$ the standard parametrization of
the rulings of the isotropic cone we have
$X(u,v)=\frac{1}{u-v}Y(v)+\frac{1}{2}Y'(v)=\frac{1}{u-v}Y(u)-\frac{1}{2}Y'(u),\
Y(v)\times Y'(v)=2iY(v)$ and with $y:=\sqrt{A}x_z$ from the TC
$y^TX(u_0,v_0)=1$ we get
$X_0=X(u_0,v_0)=\frac{Y(v_0)+iY(v_0)\times y}{y^TY(v_0)}
=\frac{Y(u_0)-iY(u_0)\times y}{y^TY(u_0)}$; replacing this into
(\ref{eq:int}) and with
$\mathcal{M}:=\pa_{u_1}y\pa_{v_1}y^T+\pa_{v_1}y\pa_{u_1}y^T
-y\pa_{u_1v_1}^2y^T$ we get
$0=Y(v_0)^TA^{-1}[y^TY(v_0)\mathcal{M}+(Y(v_0)\pm iY(v_0)\times
y)\pa_{u_1v_1}^2y^T](Y(v_0)\pm iY(v_0)\times
y)=Y(v_0)^TA^{-1}([y^TY(v_0)\mathcal{M}+Y(v_0)\pa_{u_1v_1}^2y^T]Y(v_0)
-(Y(v_0)\times y)\pa_{u_1v_1}^2y^T\\(Y(v_0)\times y)\pm
i[[y^TY(v_0)\mathcal{M}+Y(v_0)\pa_{u_1v_1}^2y^T](Y(v_0)\times
y)+(Y(v_0)\times y)\pa_{u_1v_1}^2y^TY(v_0)]),\ \forall
v_0\in\mathbb{C}$, that is two polynomials of degree $6$ in $v_0$
are identically $0$; this imposes a linear homogeneous system of
$14$ equations in $12$ variables the entries of $\mathcal{M}$ and
$\pa_{u_1v_1}^2y$ (optionally one can consider
$\mathcal{M}:=\pa_{u_1}y\pa_{v_1}y^T+\pa_{v_1}y\pa_{u_1}y^T=\mathcal{M}^T$
and we have only $9$ variables) with obvious solution
$\pa_{u_1v_1}^2y=My,\ \mathcal{M}=-M(I_3-zA),\ z\in\mathbb{C}$
discussed in \S \ref{subsec:tcsym}, so $x_z$ is ((isotropic)
singular) quadric doubly ruled by degenerate leaves and confocal
to $x_0$ (requiring a space of solutions at least $3$-dimensional
leads to over-determinate conditions on $y$).

For $x_0$ canonical (I)QWC  we have $x_0=LZ(u_0,v_0),\
L:=(\sqrt{a_1^{-1}})^{-1}e_1e_1^T+(\sqrt{a_2^{-1}})^{-1}e_2e_2^T+e_3e_3^T,\
a_1,a_2\in\mathbb{C}\setminus\{0\},\
:=(\sqrt{a_1^{-1}})^{-1}(f_1\bar f_1^T+\bar
f_1f_1^T+\frac{1}{2a_1}J_2)+e_3e_3^T,\
a_1\in\mathbb{C}\setminus\{0\}$ for QWC or $:=f_1e_3^T+\bar
f_1e_1^T-\sqrt{a_1}e_3e_2^T,\ a_1\in\mathbb{C}\setminus\{0\},\
:=J_3+\bar J_3^2$ for two IQWC; in all cases $N_0$ is a multiple
of $(L^T)^{-1}(u_0f_1+v_0\bar f_1-e_3)$ and with $y:=L^{-1}x_z$
from the TC $y^T(u_0f_1+v_0\bar f_1-e_3)=u_0v_0$ we get
$Z_0=Z(u_0,v_0)=\frac{y^T(e_3-v_0\bar
f_1)}{y^Tf_1-v_0}(f_1+v_0e_3)+v_0\bar
f_1=\frac{y^T(e_3-u_0f_1)}{y^T\bar f_1-u_0}(\bar
f_1+u_0e_3)+u_0f_1$; replacing this into (\ref{eq:int}) and with
$\mathcal{M}:=\pa_{u_1}y\pa_{v_1}y^T+\pa_{v_1}y\pa_{u_1}y^T
-y\pa_{u_1v_1}^2y^T-\pa_{u_1v_1}^2yy^T$ we get
$0=(f_1+v_0e_3)^TL^TL[(y^Tf_1-v_0)[\mathcal{M}+v_0(\bar
f_1\pa_{u_1v_1}^2y^T+\pa_{u_1v_1}^2y\bar f_1^T)]+y^T(e_3-v_0\bar
f_1)[(f_1+v_0e_3)\pa_{u_1v_1}^2y^T
+\pa_{u_1v_1}^2y(f_1+v_0e_3)^T]][y^T(e_3-v_0\bar
f_1)f_1+(y^Tf_1-v_0)(v_0\bar f_1-e_3)]=(\bar
f_1+u_0e_3)^TL^TL[(y^T\bar f_1-u_0)
[\mathcal{M}+u_0(f_1\pa_{u_1v_1}^2y^T+\pa_{u_1v_1}^2yf_1^T)]
+y^T(e_3-u_0f_1)[(\bar f_1+u_0e_3)\pa_{u_1v_1}^2y^T
+\pa_{u_1v_1}^2y(\bar f_1+u_0e_3)^T]][y^T(e_3-u_0f_1)\bar
f_1+(y^T\bar f_1-u_0)(u_0f_1-e_3)],\ \forall
u_0,v_0\in\mathbb{C}$, that is two polynomials of degree $5$
respectively in $u_0,v_0$ are identically $0$; this imposes a
linear homogeneous system of $12$ equations in the $9$ entries of
$\mathcal{M}$ and $\pa_{u_1v_1}^2y$ (optionally one can consider
$\mathcal{M}:=\pa_{u_1}y\pa_{v_1}y^T+\pa_{v_1}y\pa_{u_1}y^T$) with
obvious solution $\pa_{u_1v_1}^2y=Me_3,\ \mathcal{M}=M(f_1\bar
f_1^T+\bar f_1f_1^T-z(LL^T)^{-1}),\ z\in\mathbb{C}$ discussed in
\S \ref{subsec:tcsym}, so $x_z$ is ((isotropic) singular) quadric
doubly ruled by degenerate leaves and confocal to $x_0$ (requiring
a space of solutions at least $3$-dimensional leads to
over-determinate conditions on $y$).

When making hypothesis on the auxiliary surface $x_z$ being
(isotropic) plane or quadric one cannot use directly
(\ref{eq:int}), since it is strongly rigid (it assumes the curves
given by collapsed leaves).

If $x_z$ is a(n isotropic) plane, then we can take $e_1^Tx_z=0$
(this case is due to Bianchi \cite{B3}) or $f_1^Tx_z=0$. We shall
reproduce here Bianchi's argument, which will also apply to $x_z$
isotropic plane.

Consider
$x_0=x_0(\mathbf{x},\mathbf{y})=\mathbf{x}e_1+\mathbf{y}e_2+\mathbf{z}e_3,\
\hat
N_0:=\pa_{\mathbf{x}}\mathbf{z}e_1+\pa_{\mathbf{y}}\mathbf{z}e_2-e_3$
for $e_1^Tx_z=0$ or $x_0=\mathbf{y}f_1+\mathbf{x}\bar
f_1+\mathbf{z}e_3,\ \hat
N_0=\pa_{\mathbf{x}}\mathbf{z}f_1+\pa_{\mathbf{y}}\mathbf{z}\bar
f_1-e_3$ for $f_1^Tx_z=0$; in both cases we have
$V:=x_z-x_0=-\mathbf{x}\pa_{\mathbf{x}}x_0+w\pa_{\mathbf{y}}x_0$
and with $m=:V\times N_0+\mathbf{x}\Theta\hat N_0$ from
(\ref{eq:dissymTC}) we get $\Theta^2$ polynomial of degree at most
two in $w$. If $\Theta^2$ is a perfect square (including a-priori
the case when it does not depend of $w$), then as $w$ varies the
facets envelope two lines passing through $x_0$ and two foci in
$x_z$ (one of them possibly situated at $\infty$), so the leaves
are the lines passing through the two foci in $x_z$; if $\Theta^2$
is not a perfect square (including a-priori the case when it
depends linearly on $w$), then as $w$ varies the facets envelope a
quadratic cone centered at $x_0$ (to see this from the cone
$(t,w)\rightarrow x_0+t(c(w)-x_0),\ c(w)=c_1(w)e_2+c_2(w)e_3$ or
$=c_1(w)f_1+c_2(w)e_3$ having normal fields
$m=wv_1+v_2+\sqrt{Aw^2+2Bw+C}v_3,\ B^2-AC\neq 0,\ v_1,v_2,v_3$
linearly independent we get $c_1(w)=\frac{(x_0\times
e_3)^T(m\times\pa_w m)}{e_1^T(m\times\pa_w m)},\
c_2(w)=\frac{(e_2\times x_0)^T(m\times\pa_w m)}{e_1^T(m\times\pa_w
m)}$ or $c_1(w)=i\frac{(x_0\times e_3)^T(m\times\pa_w
m)}{f_1^T(m\times\pa_w m)},\ c_2(w)=i\frac{(f_1\times
x_0)^T(m\times\pa_w m)}{f_1^T(m\times\pa_w m)}$; in both cases
$c(w)$ is a conic), so the leaves are the tangents to a conic in
$x_z$; in both cases we have $x_z$ (isotropic) singular quadric
doubly ruled by degenerate leaves, so from \S \ref{subsec:tcsym}
$x_0$ is a quadric confocal to $x_z$ (except for the case when
$x_z$ does not admit confocal family).

Given the auxiliary surface $x_z$ being a quadric, we should get
the defining surface $x_0$ a confocal quadric (and thus by the
previous argument $x_z$ doubly ruled by degenerate leaves).

More generally given the auxiliary surface $x_z=x_z(s,w)$ one
solves the TC for $s=s(x_0,N_0,w)$; replacing this into $x_z$ we
get $V=V(x_0,N_0,w)$ and by implicit differentiation we get
$\pa_ws=-\frac{\pa_wx_z^TN_0}{\pa_sx_z^TN_0},\
ds=-\frac{V^TdN_0}{\pa_sx_z^TN_0},\
d(V+x_0)=-\frac{V^TdN_0}{\pa_sx_z^TN_0}\pa_sx_z,\
\pa_wV=\pa_wx_z-\frac{\pa_wx_z^TN_0}{\pa_sx_z^TN_0}\pa_sx_z$; the
first equation of (\ref{eq:dissymTC1}) becomes
$0=\pa_w[V^TdN_0\frac{N_0^T(\pa_sx_z\times\pa_wx_z)}{V^T(\pa_sx_z\times\pa_wx_z)}]
\wedge V^T(N_0\times
dx_0)-V^TdN_0\frac{N_0^T(\pa_sx_z\times\pa_wx_z)}{V^T(\pa_sx_z\times\pa_wx_z)}
\wedge\pa_wV^T(N_0\times
dx_0)=^{(\ref{eq:che})}\frac{N_0^T(\pa_sx_z\times\pa_wx_z)}
{V^T(\pa_sx_z\times\pa_wx_z)}N_0^T(\pa_wV\times
V)N_0^T[dN_0\times\wedge(N_0\times
dx_0)]+\pa_w\frac{N_0^T(\pa_sx_z\times\pa_wx_z)}
{V^T(\pa_sx_z\times\pa_wx_z)}V^TdN_0\wedge V^T(N_0\times
dx_0)=\pa_w\frac{N_0^T(\pa_sx_z\times\pa_wx_z)}
{V^T(\pa_sx_z\times\pa_wx_z)}V^TdN_0\wedge V^T(N_0\times dx_0)$.
If $V^TdN_0\wedge V^T(N_0\times dx_0)=0$, then from the first
equation of (\ref{eq:dissymTC}) we get $K|m|^2=0$, so we need
$\pa_w\frac{N_0^T(\pa_sx_z\times\pa_wx_z)}
{V^T(\pa_sx_z\times\pa_wx_z)}=0,\ \forall w\in\mathbb{C}$
(equivalently $\frac{N_0^TN_z}{V^TN_z}(\mod \mathrm{TC})$ does not
depend on $x_z$).

For $x_z=c_1+sc_2,\ c_j=c_j(w),\ j=1,2$ ruled with normal field
$\hat N_z=c_2\times(c_1'+sc_2')$ from the TC we get
$s=-\frac{(c_1-x_0)^TN_0}{c_2^TN_0}$ and the fraction
$\frac{N_0^T[c_2\times(c_2^TN_0c_1'-(c_1-x_0)^TN_0c_2')]}
{(c_1-x_0)^T[c_2\times(c_2^TN_0c_1'-(c_1-x_0)^TN_0c_2')]}$ should
be independent of $w$.

For $x_z=(\sqrt{A})^{-1}X(s,w)$ or $=(\sqrt{A})^{-1}X(w,s)$
canonical QC
we have $c_2(w)=(\sqrt{A})^{-1}Y(w),\\
c_1(w)=\pm\frac{1}{2}(\sqrt{A})^{-1}Y'(w)$ and we need
$F_\pm:=\frac{Y(w)^T(N_0^Tx_0\sqrt{A}N_0\pm
i(\sqrt{A})^{-1}N_0\times\sqrt{A}N_0)}
{Y(w)^T(N_0^Tx_0\sqrt{A}x_0\pm
i(\sqrt{A})^{-1}N_0\times\sqrt{A}x_0-(\sqrt{A})^{-1}N_0)}$ to be
independent of $w$, that is $F_\pm(N_0^Tx_0\sqrt{A}x_0\pm
i(\sqrt{A})^{-1}N_0\times\sqrt{A}x_0-(\sqrt{A})^{-1}N_0)=N_0^Tx_0\sqrt{A}N_0\pm
i(\sqrt{A})^{-1}N_0\times\sqrt{A}N_0$; multiplying on the left
with $N_0^T(\sqrt{A})^{-1}$ we get
$F_\pm=\frac{N_0^Tx_0}{(N_0^Tx_0)^2-N_0^TA^{-1}N_0}$ and
$(I_3-N_0N_0^T)(A^{-1}-x_0x_0^T)N_0=0$, that is the first equation
of (\ref{eq:int1}), so $x_0$ is a quadric confocal to $x_z$.

For $x_z=LZ(s,w)$ or $=LZ(w,s)$ canonical (I)QWC we have
$c_2(w)=L(\frac{e_1\mp ie_2}{\sqrt{2}}+we_3),\
c_1(w)=wL\frac{e_1\pm ie_2}{\sqrt{2}}$ and we need
$\frac{(L^{-1}N_0)^T(e_1\pm ie_2)}{(L^{-1}x_0)^T(e_1\pm
ie_2)+\frac{N_0^TL(e_1\pm ie_2)}{N_0^TLe_3}}=\frac{1}{F_\pm},\
\frac{(L^{-1}N_0)^T(\pm ie_3\times
L^TN_0-N_0^TLe_3e_3)}{(L^{-1}x_0)^T(\pm ie_3\times
L^TN_0-N_0^TLe_3e_3)-N_0^Tx_0}=\frac{1}{F_\pm},\
\frac{(L^{-1}N_0)^T[N_0^TL(e_1\mp ie_2)e_3-N_0^Tx_0(e_1\mp
ie_2)]}{(L^{-1}x_0)^T[N_0^TL(e_1\mp ie_2)e_3-N_0^Tx_0(e_1\mp
ie_2)]}=\frac{1}{F_\pm}$; from the first and last relation we get
$(L^{-1}x_0)^T(e_1\pm ie_2)=F_\pm(L^{-1}N_0)^T(e_1\pm
ie_2)-\frac{N_0^TL(e_1\pm ie_2)}{N_0^TLe_3},\
(L^{-1}x_0)^Te_3=-\frac{N_0^Tx_0}{N_0^TLe_3}+F_\mp(L^{-1}N_0)^Te_3
+(F_\pm-F_\mp)N_0^Tx_0\frac{(L^{-1}N_0)^T(e_1\pm
ie_2)}{N_0^TL(e_1\pm ie_2)}$; replacing this into the second
relation we get $(L^{-1}x_0)^T(e_3\times
L^TN_0)=F_\pm(L^{-1}N_0)^T(e_3\times L^TN_0)\pm
i(F_\pm-F_\mp)N_0^Tx_0\frac{(L^{-1}N_0)^T(e_1\mp
ie_2)}{N_0^TL(e_1\mp ie_2)}N_0^TLe_3$; applying compatibility
condition to the last two relations we get
$0=(F_+-F_-)(L^{-1}N_0)^T(e_3N_0^TL-2N_0^Tx_0I_3)(f_1\bar
f_1^T+\bar f_1f_1^T)L^TN_0 =(F_+-F_-)(L^{-1}N_0)^T(e_3\times
L^TN_0)N_0^T(L(f_1\bar f_1^T+\bar
f_1f_1^T)L^T+x_0(Le_3)^T+Le_3x_0^T)N_0$.

If $F_+\neq F_-$, then the conditions imposed on $x_0$ are
over-determined; otherwise we get $(I_3-N_0N_0^T)(L(f_1\bar
f_1^T+\bar f_1f_1^T)L^T+x_0(Le_3)^T+Le_3x_0^T)N_0=0$, that is the
second equation of (\ref{eq:int1}), so $x_0$ is a quadric confocal
to $x_z$.

It remains to show that these are all the cases of B
transformation with defining surface.

Multiplying (\ref{eq:int}) on the left with $N_z^T$ we get

$\frac{N_z^TN_0}{N_z^TV}=-\frac{\pa_{u_1v_1}^2x_z^TN_0}{2\pa_{u_1}x_z^TN_0\pa_{v_1}x_z^TN_0}$;
differentiating this with respect to $v_1$ and keeping account of
$\pa_{v_1}u_1=-\frac{\pa_{v_1}x_z^TN_0}{\pa_{u_1}x_z^TN_0}$ we get
\begin{eqnarray}\label{eq:pauv}
\pa_{u_1v_1v_1}^3x_z^TN_0\pa_{u_1}x_z^TN_0-\pa_{u_1u_1v_1}^3x_z^TN_0\pa_{v_1}x_z^TN_0=
\nonumber\\
\frac{\pa_{u_1v_1}^2x_z^TN_0}{\pa_{u_1}x_z^TN_0\pa_{v_1}x_z^TN_0}[
(\pa_{u_1}x_z^TN_0)^2\pa_{v_1}^2x_z^TN_0-(\pa_{v_1}x_z^TN_0)^2\pa_{u_1}^2x_z^TN_0]
\end{eqnarray}

Differentiating (\ref{eq:int}) with respect to $v_1$ we get
$0=(I_3-N_0N_0^T)[\pa_{u_1}x_z^TN_0(\pa_{u_1}x_z\pa_{v_1}^2x_z^T+\pa_{v_1}^2x_z\pa_{u_1}x_z^T)
-\pa_{v_1}x_z^TN_0(\pa_{v_1}x_z\pa_{u_1}^2x_z^T+\pa_{u_1}^2x_z\pa_{v_1}x_z^T)]N_0-
V[\pa_{u_1v_1v_1}^3x_z^TN_0\pa_{u_1}x_z^TN_0-\pa_{u_1u_1v_1}^3x_z^TN_0\pa_{v_1}x_z^TN_0]=
^{(\ref{eq:pauv})}(I_3-N_0N_0^T)[\pa_{u_1}x_z^TN_0(\pa_{u_1}x_z\pa_{v_1}^2x_z^T+\pa_{v_1}^2x_z\pa_{u_1}x_z^T)
-\pa_{v_1}x_z^TN_0(\pa_{v_1}x_z\pa_{u_1}^2x_z^T+\pa_{u_1}^2x_z\pa_{v_1}x_z^T)]N_0\\-
V\frac{\pa_{u_1v_1}^2x_z^TN_0}{\pa_{u_1}x_z^TN_0\pa_{v_1}x_z^TN_0}[
(\pa_{u_1}x_z^TN_0)^2\pa_{v_1}^2x_z^TN_0-(\pa_{v_1}x_z^TN_0)^2\pa_{u_1}^2x_z^TN_0]$;
replacing $V\pa_{u_1v_1}^2x_z^TN_0$ from (\ref{eq:int}) we get
\begin{eqnarray}\label{eq:fina}
V^TN_0=0\Rightarrow[\frac{\pa_{u_1}^2x_z\times\pa_{u_1}x_z}{(\pa_{u_1}x_z^TN_0)^3}-
\frac{\pa_{v_1}^2x_z\times\pa_{v_1}x_z}{(\pa_{v_1}x_z^TN_0)^3}]\times
N_0=0.
\end{eqnarray}
If
$\frac{\pa_{u_1}^2x_z\times\pa_{u_1}x_z}{(\pa_{u_1}x_z^TN_0)^3}-
\frac{\pa_{v_1}^2x_z\times\pa_{v_1}x_z}{(\pa_{v_1}x_z^TN_0)^3}\neq
0$, then along the curve of tangency of the tangent cone of $x_0$
from $x_z(u_1,v_1)$ $N_0$ spans a line or a plane; since the
tangent cone cannot be a cylinder $x_0$ must be developable.

Thus
$\frac{\pa_{u_1}^2x_z\times\pa_{u_1}x_z}{(\pa_{u_1}x_z^TN_0)^3}-
\frac{\pa_{v_1}^2x_z\times\pa_{v_1}x_z}{(\pa_{v_1}x_z^TN_0)^3}=0$;
since along the curve of tangency of the tangent cone of $x_0$
from $x_z(u_1,v_1)$ for $u_1,v_1$ fixed
$\pa_{u_1}x_z^TN_0\pa_{v_1}x_z^TN_0$ is constant, the ratio
$\pa_{u_1}x_z^TN_0/\pa_{v_1}x_z^TN_0$ cannot be constant, so
$\pa_{u_1}^2x_z\times\pa_{u_1}x_z=\pa_{v_1}^2x_z\times\pa_{v_1}x_z=0$
and $x_z$ is doubly ruled by degenerate leaves.

\end{document}